\documentclass[]{article}

\usepackage{hyperref}

\renewcommand{\v}[1]{\ensuremath{\mathbf{#1}}}

\newcommand{\mc}{\mathcal}

\def\st{{\rm s.t.}}

\usepackage[authoryear]{natbib}
\usepackage{algorithmic}
\usepackage{algorithm}
\usepackage{booktabs}
\usepackage{caption,subcaption}
\usepackage{comment}
\usepackage[usenames,dvipsnames]{color}
\usepackage{fullpage,graphicx,amssymb,amsmath}
\usepackage{enumerate}
\usepackage{empheq}
\usepackage[authoryear]{natbib}
\usepackage[normalem]{ulem}
\usepackage{url}
\usepackage{siunitx} 

\usepackage{multirow}
\usepackage{enumitem}

\renewcommand{\v}[1]{\ensuremath{\mathbf{#1}}}
\def\st{{\rm s.t.}}

\newcommand{\Cvar}{\mathbb{CVAR}} 
\newcommand{\Var}{\mathbb{VAR}} 

\newcommand{\Cov}{{\rm Cov}} 
\newcommand{\E}{\mathbb{E}} 
\renewcommand{\P}{\mathbb{P}} 
\renewcommand{\Re}{\mathbb{R}} 
\newcommand{\T}{\mathcal{T}} 

\newcommand{\vell}{\boldsymbol{\ell}}

\newcommand{\vxi}{\boldsymbol{\xi}}

\newcommand{\vSigma}{\boldsymbol{\Sigma}}

\newcommand{\term}{x^{\textrm{Term}}}
\newcommand{\spot}{y^{\textrm{Spot}}}
\newcommand{\prodCost}{C^{\textrm{Prod}}}
\newcommand{\transCost}{C^{\textrm{Trans}}}
\newcommand{\prodVar}{u^{\textrm{Prod}}}
\newcommand{\transVar}{u^{\textrm{Trans}}}

\newcommand{\cmt}[1]{\textcolor{black}{#1}}

\title{Data-Driven Distributionally Robust Optimization for\\ Long-Term Contract vs. Spot Allocation Decisions:\\ Application to Electricity Markets}
\author{Dimitri J. Papageorgiou \\
{\small ExxonMobil Technology and Engineering Company} \\
{\small 1545 Route 22 East, Annandale, NJ 08801 USA} \\
{\small dimitri.j.papageorgiou@exxonmobil.com} \\
\\
}

\begin{document}

\maketitle

\begin{abstract}
There are numerous industrial settings in which a decision maker must decide whether to enter into long-term contracts to guarantee price (and hence cash flow) stability or to participate in more volatile spot markets. 
In this paper, we investigate a data-driven distributionally robust optimization (DRO) approach aimed at  balancing this tradeoff.  Unlike traditional risk-neutral stochastic optimization models that assume the underlying probability distribution generating the data is known, DRO models assume the distribution belongs to a family of possible distributions, thus providing a degree of immunization against unseen and potential worst-case outcomes.
We compare and contrast the performance of a risk-neutral model, conditional value-at-risk formulation, and a Wasserstein distributionally robust model to demonstrate the potential benefits of a DRO approach for an ``elasticity-aware'' price-taking decision maker. 

\vspace{0.2in}

\parindent=0cm \textit{Keywords}

\parindent=0.5cm Conditional value-at-risk, Data-driven distributionally robust optimization, Price elasticity, Wasserstein metric.
\end{abstract}

\section*{Nomenclature} \label{sec:nomenclature}
\begin{tabular}{ll}
\toprule
& \textbf{Definition} \\
\hline
\multicolumn{2}{l}{\textbf{Indices and sets}} \\
	$t \in \T$ 	& set of time periods: $\T = \{1,\dots,T\}$ \\	
	$m \in \mc{M}$ 	& set of market locations \\
	$c \in \mc{C}_m$	& set of contract associated with market location $m$ \\
	$i \in \mc{I}$ & set of supply steps  \\
	$\mc{K}_m$ & set of steps defining the elasticity curve (step function) in market $m$ \\
	$\mc{P}$ 	& set (or family) of distributions containing the true data-generating mechanism \\
	$s \in \mc{S}$ 	& finite set of scenarios \\
\\
\multicolumn{2}{l}{\textbf{Parameters}} \\
	$\alpha \in (0,1)$ & probability level used in conditional value-at-risk of Formulation~\eqref{model:cvar}\\ 
	$\gamma \in (0,1)$ & probability level used in conditional value-at-risk for ``reward-to-risk'' $\rho_{\gamma}(\v{y}^{\textrm{Spot}})$ metric \\ 
	$\lambda \in [0,1]$ & user-defined parameter governing the weight given to the expected value\\ 
	$\epsilon \in \Re_+$ & Wasserstein radius defines the permissible deviation from the empirical data \\ 
	$\pi_s \in (0,1]$ & probability of scenario $s$ \\ 
	$P_{mkts}$ & spot market price in market $m$  in $(t,s)$ (shorthand for ``in time period $t \in \mc{T}$ in scenario $s \in \mc{S}$.'') \\
	$\transCost_{m}$ & per-unit transportation cost to market location $m$ \\
	$\prodCost_i$ & production cost of  for each supply step $i \in \mc{I}$ \\
	$L_t$ ($U_t$) & minimum (maximum) production limits in time period $t$ \\
	$U_{i}^{\textrm{Prod}}$ & maximum production available at supply step $i \in \mc{I}$ \\
	$W_{mct}$ & long-term ``wholesale'' price associated with contract $c \in \mc{C}_m$ in market $m$ in time period $t \in \mc{T}$ \\
	$Y_{mkts}^{\textrm{Spot}}$ & maximum spot volume that can be allocated to spot market elasticity step $k$ in $(t,s)$ \\
\\
\multicolumn{2}{l}{\textbf{Continuous Decision variables}} \\
	$\prodVar_{its}$ & production amount at supply step $i \in \mc{I}$ in $(t,s)$ \\
	$\transVar_{mts}$ & amount transported to market location $m \in \mc{M}$ in $(t,s)$ \\
	$x_{mc}^{\min}$ & long-term contract volume allocation to market-contract pair $(m,c)$ \\
	$\term_{mcts}$ & production amount to allocate to long-term contract $c \in \mc{C}_m$ in $(t,s)$ \\
	$\spot_{mkts}$ & production amount to allocate to spot market in location $m$ on elasticity curve step $k$ in $(t,s)$ \\
\\
\multicolumn{2}{l}{\textbf{Functions}} \\
	$\zeta(\v{y}^{\textrm{Spot}})$ & Expected risk-neutral profit given a spot allocation decision $\v{y}^{\textrm{Spot}}$ \\
	$\chi_{\gamma}(\v{y}^{\textrm{Spot}})$ & Expected $(1-\gamma)$-tail profit ($\Cvar_{\gamma}$) given a spot allocation decision $\v{y}^{\textrm{Spot}}$ and $\gamma \in (0,1)$ \\
	$\zeta^{\textrm{riskfree}}$ & Expected risk-free profit. Note that $\zeta^{\textrm{riskfree}}=\zeta(\v{0})=\chi_{\gamma}(\v{0}) \quad \forall \gamma \in (0,1)$ \\
	$\Delta\zeta(\v{y}^{\textrm{Spot}})$ & Expected increase in the risk-neutral profit relative to the risk-free profit \\
		& given a spot allocation decision $\v{y}^{\textrm{Spot}}$: $\Delta\zeta(\v{y}^{\textrm{Spot}}) =\zeta(\v{y}^{\textrm{Spot}})-\zeta^{\textrm{riskfree}}$ \\
	$\Delta\chi_{\gamma}(\v{y}^{\textrm{Spot}})$ & Absolute value of the expected decrease in the $(1-\gamma)$-tail profit relative to the risk-free profit \\
		& given a spot allocation decision $\v{y}^{\textrm{Spot}}$: $\Delta\chi_{\gamma}(\v{y}^{\textrm{Spot}}) =|\chi_{\gamma}(\v{y}^{\textrm{Spot}})-\zeta^{\textrm{riskfree}}|$ \\
	$\rho_{\gamma}(\v{y}^{\textrm{Spot}})$ & Change in expected profit per change in risk (expected $(1-\gamma)$-tail profit)  \\
		& as a function of the spot allocation $\v{y}^{\textrm{Spot}}$: $\Delta\zeta(\v{y}^{\textrm{Spot}}) / \Delta \chi_{\gamma}(\v{y}^{\textrm{Spot}})$ \\
\bottomrule
\end{tabular}

\section{Introduction} \label{sec:intro}

In numerous energy markets, including electricity and natural gas, energy providers/sellers face the dilemma of deciding how to allocate their supply across various markets and over time \cmt{\citep{kirschen2018fundamentals,gabriel2012complementarity}}.  Throughout this paper, for convenience, we will use the motivating example of a generation company (Genco) who sells electricity. There are typically at least two major types of markets where Gencos can sell electricity: the spot market and the forward market for longer-term bilateral contracts \citep{kirschen2018fundamentals}.  Here, the spot market refers to a public financial market where electricity is traded on a daily, hourly, or subhourly basis. Gencos deliver electricity immediately and buyers pay for it ``on the spot.''  Such markets can be highly volatile as supply and demand variability can cause the market-clearing price to fluctuate dramatically over time \cmt{\citep{mayer2018electricity,shahidehpour2017restructured}}.  In contrast, forward markets allow for buyers and sellers to enter longer-term bilateral contracts to reduce price variability over a time span of interest. Among other benefits, forward markets allow market participants to hedge against uncertainty by offering price predictability, which in turn allows the Genco to better plan its cash flows and potentially secure more favorable financing from lenders. A power purchase agreement (PPA) is one such example of a bilateral contract where a Genco enters a long-term agreement to provide (a typically fixed amount of) energy at a fixed price over many time periods, e.g., one year. 
This paper investigates the problem of deciding how to allocate supply within these two types of markets.

\subsection{Electricity applications investigating long-term vs. spot tradeoffs} \label{sec:electricity_applications}

Given the importance of balancing risk exposure and expected profits in the power sector, many researchers have studied how power producers should simultaneously optimize contractual involvement, spot allocation, and generation planning. Within the electricity sector, this joint problem is often referred to as ``power portfolio optimization'' and ``integrated risk management'' \citep{lorca2014power}. Although there are a host of contracts available, in this work we focus exclusively on fixed-price forward contracts; we do not consider other options and derivatives. A forward contract is an agreement to buy or sell a fixed amount of electricity at a given price over a fixed time horizon. 

The importance of forward contracts within the electricity industry has been known for decades \citep{kaye1990forward}. We therefore attempt here to highlight key papers that have employed a rigorous approach to address power portfolio optimization under uncertainty.
Early papers employed more traditional risk-neutral stochastic programming techniques as a means to transcend a deterministic Markowitz mean-variance mindset \citep{kwon2006stochastic,sen2006stochastic}.
Risk-averse models soon followed to manage downside exposure.
Noteworthy papers that employ the conditional value-at-risk ($\Cvar$) metric include \cite{conejo2008optimal}, \cite{street2009bidding}, \cite{pineda2012managing}, \cite{yau2011financial}, \cite{lorca2014power}. 
For example, \citet{street2009bidding} investigate bidding strategies for risk-averse Gencos in a long-term forward contract auction.
\citet{fanzeres2014contracting} examine contracting strategies for renewable generators using an interesting hybrid stochastic/robust optimization approach.

\begin{figure}[h!] 
\centering
\includegraphics[width=16cm]{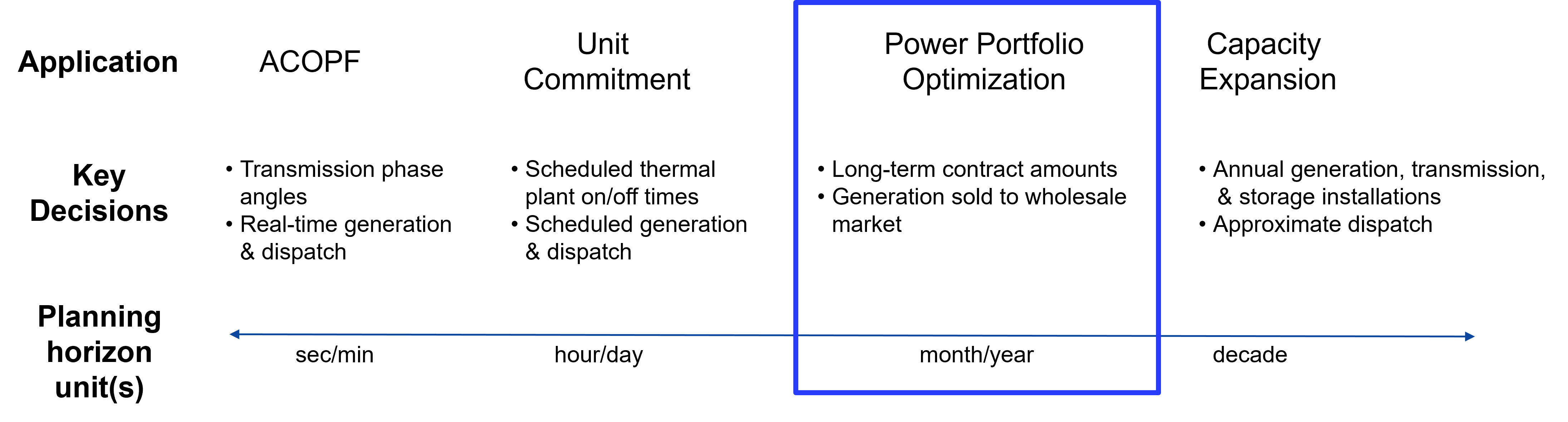}
\caption{The power portfolio optimization problem in relation to other prominent electricity planning and scheduling problems.}
\label{fig:Power_portfolio_opt_Simple}
\end{figure} 

Our investigation most closely aligns with the study conducted by \citet{lorca2014power},
who pose a risk-averse stochastic linear program for a ``medium-term'' planning horizon.  Consistent with the planning horizons shown in Figure~\ref{fig:Power_portfolio_opt_Simple}, the qualifier ``medium-term'' is used to distinguish the problem from short-term scheduling problems over minutes, hours, or days \cmt{\citep{bienstock2020mathematical,xavier2021learning}}, and from long-term applications that may involve capital-intensive generation investment decisions \cmt{\citep{lara2018deterministic,guo2022generation}}. We also consider a medium-term planning horizon of one year.

\cmt{Although not directly related to this work, several research groups have considered many facets of chemical process systems interacting with electricity markets. Noteworthy exemplars include demand respond contracts of energy intensive processes, e.g., air separation units \citep{zhang2016enterprise}; dynamic operation of energy intensive chemical processes \citep{cao2016optimal,pattison2016optimal}; scheduling of buildings/HVAC systems \citep{risbeck2017mixed,zavala2013real}; integrated energy systems and ancillary services \citep{dowling2018economic}; distributed chemical manufacturing, e.g., ammonia \citep{palys2018exploring,shao2019space}; interactions with the electric grid, e.g., bidding \citep{dowling2017multi}; scheduling of co-generation systems participating in energy markets \citep{mitra2013optimal}; dynamic operation of carbon capture systems for power generation \citep{bhattacharyya2011steady,kelley2018demand}. Many of these applications could benefit from optimization under uncertainty, assuming there is sufficient data availability and fidelity.}

\subsection{Why pursue distributionally robust optimization?}

Over the past decade, distributionally robust optimization (DRO) has emerged as a powerful tool within the operations research and statistical learning communities, while also garnering attention within the process systems engineering community (see, e.g., \cite{gao2019data,liu2021multistage,shang2018distributionally}). 
\citet{rahimian2019droreview} survey fundamental DRO concepts and applications, in addition to relating it with robust optimization, risk aversion, chance-constrained optimization, and function regularization. Loosely speaking, DRO is well-suited to address data-driven optimization problems as it puts faith in the empirical data, but not too much. 

\citet{esfahani2018dro} motivate DRO quite nicely. A traditional stochastic program attempts to solve the problem $\min_{\v{x} \in \mc{X}} \E^{\P}[h(\v{x},\vxi)]$, where the loss function $h:\Re^n \times \Re^m$ depends on both the decision vector $\v{x} \in \Re^n$ and the random vector $\vxi \in \Re^m$ governed by the distribution $\P$. Unfortunately, as many practitioners have discovered, the true distribution $\P$ is rarely known precisely and must be inferred from data, physics, expert knowledge, and more.  Optimizing a traditional stochastic program can then lead to solutions that are ``overfitted to the data.'' Second, computing the expectation in a stochastic program for a fixed decision $\v{x}$ can be computationally challenging in its own right as it may involve the evaluation of a multivariate integral.  

To combat these challenges, DRO attempts to hedge the expected loss against a family $\mc{P}$ of distributions that include the true data-generating mechanism with high confidence \citep{chen2018robust}. Mathematically, DRO minimizes the expected loss over the worst-case distribution $\mathbb{Q} \in \mc{P}$ by solving 
\begin{equation} \label{model:generic_DRO}
\min_{\v{x} \in \mc{X}} \sup_{\mathbb{Q} \in \mc{P}} \E^\mathbb {Q} \big[ h(\v{x},\vxi) \big].
\end{equation}
The family $\mc{P}$ is also known as an \textit{ambiguity set} and minimizing the inner ``sup'' expectation is sometimes termed an \textit{ambiguity-averse} (as opposed to ``risk-averse'') problem, which is why DRO is also known as \textit{ambiguous stochastic optimization} \citep{rahimian2019droreview}. DRO ``bridges
the gap between data and decision-making -- statistics and optimization frameworks -- 
to protect the decision-maker from the ambiguity in the underlying probability distribution'' \citep{rahimian2019droreview}.

\begin{table}
\caption{Pros and cons of moment- and statistical-based ambiguity sets}
\label{table:continuous_vs_discrete_time}
\begin{center}
    \begin{tabular}{ l p{6.5cm} p{6.5cm} }
    \toprule
    \textbf{} & \textbf{Moment-based ambiguity sets} & \textbf{Statistical-based ambiguity sets} \\\hline
    \multirow{3}{*}{Advantages}
    &
    \vspace{-1.0ex}
    \begin{enumerate}[noitemsep, topsep=0pt, leftmargin=*]
    \item Give rise to tractable formulations
    \item Some claim that the assumptions are easier to communicate 
    \end{enumerate}
    &
    \vspace{-1.0ex}
    \begin{enumerate}[noitemsep, topsep=0pt, leftmargin=*]
    \item Give rise to tractable formulations
    \item Data-driven
    \end{enumerate}
    \\[-2.0ex]
    \hline
    \multirow{5}{*}{Disadvantages}
    &
    \vspace{-1.0ex}
    \begin{enumerate}[noitemsep, topsep=0pt, leftmargin=*]
    \item Why should we expect to know the moment conditions, but nothing else?
    \item The resulting worst-case distributions sometimes yield overly conservative decisions
    \end{enumerate}
    &
    \vspace{-1.0ex}
    \begin{enumerate}[noitemsep, topsep=0pt, leftmargin=*]
    \item Difficult to communicate assumptions?
    \item Deciding which statistical metric to employ may be challenging
    \end{enumerate}
    \\[-2.0ex]
    \bottomrule
    \end{tabular}
\end{center}
\end{table}

Existing DRO approaches can be divided into two broad classes -- moment-based and statistical distance-based -- according to the way in which the ambiguity set $\mc{P}$ is constructed.
Table~\ref{table:continuous_vs_discrete_time} captures the main purported advantages and disadvantages of the types of ambiguity sets. 
Moment-based ambiguity sets postulate that the empirical data must belong to a distribution that satisfies certain moment (e.g., mean and variance) constraints \citep{delage2010distributionally}.  While such approaches often give rise to tractable formulations, they sometimes produce overly conservative solutions \citep{wang2016likelihood} and may not enjoy favorable asymptotic consistency or finite sample guarantees \citep{hanasusanto2018conic}.
On the other hand, statistical distance-based ambiguity sets require distributions that are stastically ``close'' to the empirical distribution.  Popular choices of distance metrics include Kullback-Leibler divergence, $\phi$-divergence, the Prokhorov metric, total variation, and more \citep{rahimian2019droreview}. Due to several shortcomings in the aforementioned metrics \citep{gao2016distributionally}, the Wasserstein distance metric has garnered considerable attention in the past decade, within both the machine learning and optimization communities.  It possesses favorable statistical guarantees, while also leading to tractable optimization formulations \citep{gao2016distributionally,esfahani2018dro,rahimian2019droreview}. It is for these reasons that we pursue the Wasserstein metric in this study.

\subsection{Contributions} \label{sec:contributions}

The contributions of this paper are: 
\begin{enumerate}
\item In contrast to the prevailing simplistic price-taker models commonly found in the literature, we consider an ``elasticity-aware'' decision-maker. Consequently, the supplier can estimate the price elasticity due to her own supply to each spot market and behaves accordingly so as not to oversaturate a particular market and ultimately overly depress prices. 
\item We present a data-driven DRO approach exploiting Wasserstein ambiguity sets and contrast it against a standard conditional value-at-risk approach. As data-driven techniques that bridge data science, machine learning, and optimization hold significant promise, we believe that this investigation is valuable for the process systems engineering community. 
\item We provide numerical evidence that our risk-averse models are tractable for an interesting application in the real-time PJM electricity market. We also explore the tradeoff between supply allocation to long-term contracts vs. the spot market as a function of the decision maker's risk aversion.  
\end{enumerate}

The remainder of this paper is organized as follows.
In addition to describing the power portfolio optimization problem aimed at optimally allocating supply across markets and customers, Section~\ref{sec:problem_statement} introduces three different two-stage stochastic programming formulations to tackle the problem: a risk-neutral stochastic program, a risk-averse model using conditional value-at-risk, and a DRO model over a Wasserstein ball.
Section~\ref{sec:numerical_results} presents numerical results based on two case studies involving the PJM electricity market. 
Conclusions and future research directions are offered in Section~\ref{sec:conclusions}.

\section{Problem Statement and Formulations} \label{sec:problem_statement}

\begin{figure}[h!] 
\centering
\includegraphics[width=16cm]{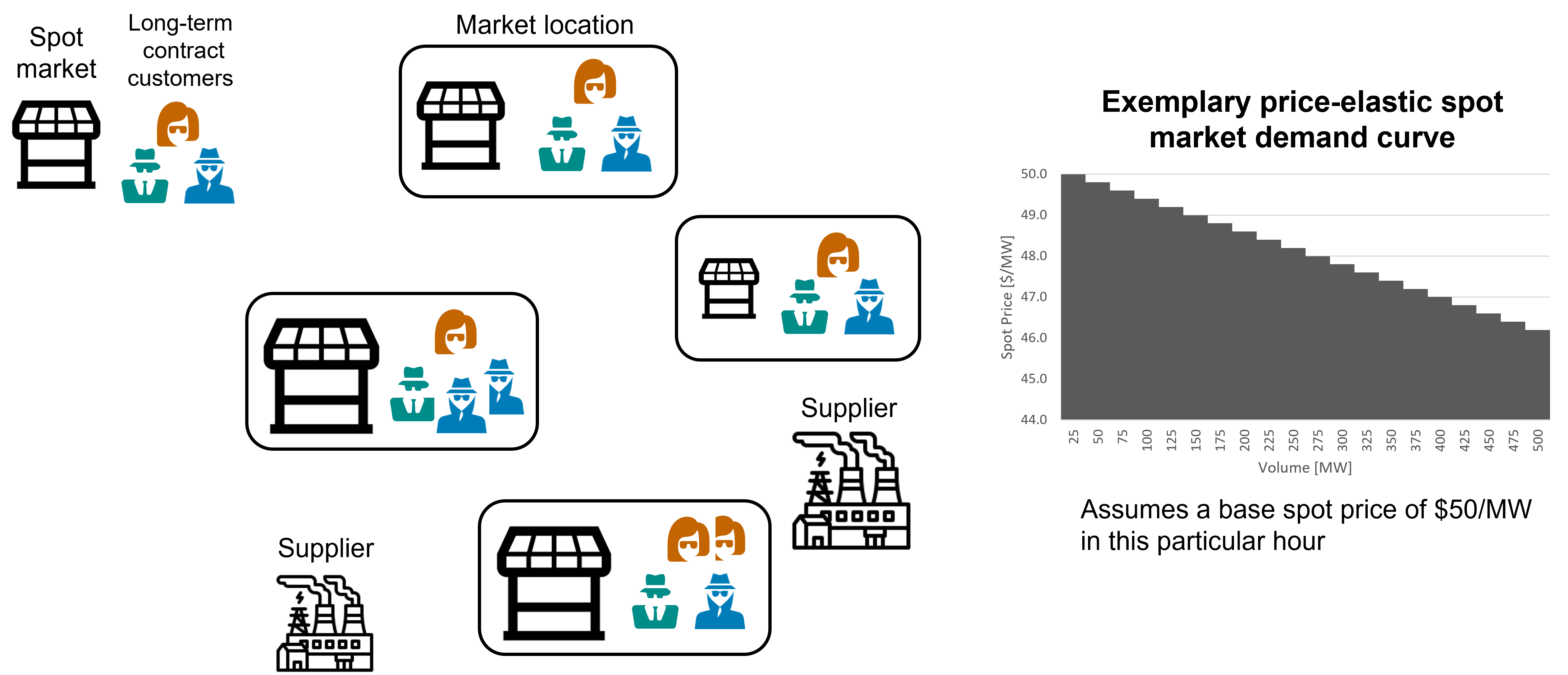}
\caption{Illustration of key ingredients in the power portfolio optimization problem.}
\label{fig:Power_portfolio_opt_problem_statement}
\end{figure}

Suppose that a key decision maker's objective is to maximize profit by selling a commodity (e.g., electricity) in a set $\mc{M}$ of market locations over a fixed planning horizon. 
As shown in Figure~\ref{fig:Power_portfolio_opt_problem_statement}, within each market location $m \in \mc{M}$, she has two options available: (1) enter into long-term fixed-price forward bilateral contracts with customers to ensure price predictability, or (2) sell to one or more spot markets at a potentially volatile market-clearing price.
More formally, she must choose a long-term contract amount $x_{mc}^{\min} \in \Re_+$ (a non-negative scalar) for each market location $m \in \mc{M}$ and each long-term contract $c \in \mc{C}_m$, which holds for the entire planning horizon (e.g., one month).  The contracted amount $x_{mc}^{\min}$ determines the minimum and maximum amount of a particular product/commodity that must be sold at a fixed price to the associated contract holder in each time period (e.g., hour) $t \in \mc{T}$.  
Specifically, she can sell up to $x_{mc}^{\min} + X_{mc}^+$ to long-term contracts where $X_{mc}^+ \in \Re_+$. 
Any remaining production in that time period can be sold to one or more spot markets.
Let $W_{mct}$ be the long-term price (sometimes called a ``wholesale'' price) associated with contract $c \in \mc{C}_m$ in market $m$ in time period $t \in \mc{T}$. Oftentimes, the long-term contract price is not a function of time and can therefore be written simply as $W_{mc}$, but we will keep the subscript $t$ for the more general setting in which the price is known to vary with time.

Meanwhile, we assume that spot market prices in each market location $m$ are unknown when deciding long-term contract amounts $x_{mc}^{\min}$. In a data-driven setting, we assume that we have access to a finite set $\mc{S}$ of scenarios, e.g., a set of historical spot price time series.  Associated with each scenario $s \in \mc{S}$ is a probability $\pi_{s} \in (0,1]$ and spot price curve for every market $m$ and time period $t \in \mc{T}$. Naturally, $\sum_{s \in \mc{S}} \pi_s = 1$. 
Unlike the majority of price-taker models in the literature, we assume that the decision maker is an ``elasticity-aware'' price-taking supplier, i.e., the supplier can estimate the price elasticity due to her own supply to each spot market $m$.  We do not assume price setter behavior, however, in which suppliers could exert market power and could therefore act as Nash-Cournot players \cmt{\citep{ralph2006EPECs,kazempour2010strategic}}. 
In our simplified setting and our numerical experiments, we assume that there is no cross-market price elasticity. \cmt{That is, we assume that the stochastic spot price in market location $m_1$ is independent of the price (and quantity sold) in market location $m_2$ for all $m_1 \neq m_2 \in \mc{M}$.} This simplification \cmt{corresponds to the setting in which all markets are independent and} allows us to represent the spot market price in each location via a descending ``staircase'' structure (shown on the right of Figure~\ref{fig:Power_portfolio_opt_problem_statement}) defined by a set $\mc{K}_m$ of steps.  The height of step $k$ in time period $t$ in scenario $s$ is denoted by $P_{mkts}$ and satisfies $P_{m1ts} > P_{m2ts} > \cdots > P_{mKts}$ 
(with $K=|\mc{K}_m|$), and has width $Y_{mkts}^{\textrm{Spot}} > 0$. 
This spot price structure implies that at most $Y_{m1ts}^{\textrm{Spot}}$ units can be sold to spot at price $P_{m1ts}$ before the price decreases to $P_{m2ts}$ and so forth.
We later describe an extension to handle ``separable'' cross-market price elasticity.


We assume that the decision maker also has production and transportation costs to consider.
Analogous to the aforementioned demand curves, we assume a standard supply curve represented by an increasing ``staircase'' structure.  That is, the supplier can produce up to $U_{i}^{\textrm{Prod}}$ units at a production cost of $\prodCost_i$ for each supply step $i \in \mc{I}$.
Let $L_t$ and $U_t$ be the known minimum and maximum production limits in time period $t$.
Let $\transCost_{m}$ denote the per-unit transportation cost to market location $m$.
Our basic model implicitly assumes that all supply is co-located, e.g., at a centralized location; this assumption could easily be relaxed.

\cmt{We do not model transmission constraints that could impact the amount of supply that could be allocated to a particular market. Such constraints could diminish the amount of supply that is transmitted to a market in a high spot price scenario. Moreover, such constraints could be time-dependent as it is well known that in warmer temperatures, transmission lines generally carry less power, which can contribute to tight grid conditions \cite{ERCOT2020}. Finally, such constraints could be partly responsible for price spikes as they lead to an imbalance in energy being available in the right place at the right time.} 

As for the decision variables, let $\term_{mcts}$ and $\sum_{k \in \mc{K}_m} \spot_{mkts}$ denote the amount of production to allocate to long-term contract $c \in \mc{C}_m$ and to the spot market in location $m$, respectively, in time period $t$ in scenario $s$. As stated above, $x_{mc}^{\min}$ denotes the long-term contract volume allocation to market-contract pair $(m,c)$. Finally, $\prodVar_{its}$ and $\transVar_{mts}$ denote the production amount at supply step $i \in \mc{I}$ and the amount transported to market $m$ in time period $t$ in scenario $s$.

\subsection{Risk-neutral stochastic programming formulation}

We first provide a basic risk-neutral stochastic programming formulation.  In all of our models, we assume that a finite set $\mc{S}$ of scenarios is available to us.
Assuming no cross-market price elasticity, a potential scenario-based risk-neutral stochastic mixed-integer linear program (MILP) has the following form:
\begin{subequations} \label{model:risk_neutral_basic}
 \begin{align}
 \max_{\substack{\v{x}^{\min},\\\v{u},\v{x},\v{y},\v{z}}}~~& \sum_{s \in \mc{S}} \pi_{s} z_s \quad (= \textrm{Expected Profit}) \label{objfnc:risk_neutral} \\
 \st~~
 	& z_s = \sum_{t \in \mc{T}} \Bigg[ \sum_{m \in \mc{M}} \sum_{c \in \mc{C}_m} W_{mct} \term_{mcts} + \sum_{k \in \mc{K}_m} P_{mkts} \spot_{mkts}  \notag \\
	& \quad - \sum_{i \in \mc{I}} \prodCost_i \prodVar_{its} - \sum_{m \in \mc{M}} \transCost_{m} \transVar_{mts} \Bigg] \quad \forall s \in \mc{S} \label{eq:profit_in_scenario_s} \\
	& x_{mc}^{\min} \leq \term_{mcts} \leq x_{mc}^{\min} + X_{mc}^+ \quad \forall m,c \in \mc{C}_m,t \in \mc{T}, s \in \mc{S} \label{eq:risk_neutral_wholesale_allocation_contract_requirement} \\	
	& \sum_{i \in \mc{I}} \prodVar_{its} = \sum_{m \in \mc{M}} \sum_{c \in \mc{C}_m} \term_{mcts} + \sum_{k \in \mc{K}_m} \spot_{mkts} \quad \forall t \in \mc{T}, s \in \mc{S} \label{eq:supply_demand_balance}\\	
	& \sum_{i \in \mc{I}} \prodVar_{its} = \sum_{m \in \mc{M}} \transVar_{mts}  \qquad \forall t \in \mc{T}, s \in \mc{S} \label{eq:supply_transportation_balance} \\	
	& L_t \leq \sum_{i \in \mc{I}} \prodVar_{its} \leq U_t \qquad \forall t \in \mc{T},s \in \mc{S} \label{eq:supply limits}  \\
	& (\v{u},\v{x}^{\min},\v{x}^{\textrm{Term}},\v{y}^{\textrm{Spot}},\v{z}) \in \mc{X}  \qquad  \label{eq:side_constraints}  \\
	& \prodVar_{its} \in [0,U_{i}^{\textrm{Prod}}] \quad \forall i \in \mc{I}, t \in \mc{T}, s \in \mc{S} \\
	& \transVar_{mts} \in [0,U_t] \quad \forall m \in \mc{M}, t \in \mc{T}, s \in \mc{S} \\	
	& x_{mc}^{\min} \in [0,X_{mc}^{\max}]  \quad \forall m \in \mc{M},c \in \mc{C}_m \\	
	& \term_{mcts} \in [0,X_{mc}^{\max}]  \quad \forall m \in \mc{M},c \in \mc{C}_m, t \in \mc{T}, s \in \mc{S} \\	
	& \spot_{mkts} \in [0,Y_{mkts}^{\textrm{Spot}}] \quad \forall m \in \mc{M}, k \in \mc{K}_m, t \in \mc{T},s \in \mc{S} \label{eq:spot_variable_bounds} \\
	& z_{s} \in \Re \qquad \forall s \in \mc{S} \label{eq:profit_in_scenario_s_variable_bounds}
\end{align}
\end{subequations}
The objective function~\eqref{objfnc:risk_neutral} is the expected or ``sample average'' profit over a finite set of scenarios. The profit $z_s$ in scenario $s$ in equation~\eqref{eq:profit_in_scenario_s} includes two terms: the two positive terms account for revenues from the long-term and spot markets, while the two negative terms denote the cost to produce and transport volumes to markets. As a reminder, the only uncertain parameter is the spot price $P_{mkts}$, rendering Formulation~\eqref{model:risk_neutral_basic} a stochastic MILP with objective function uncertainty only. 
Constraints~\eqref{eq:risk_neutral_wholesale_allocation_contract_requirement} ensure that the amount of production allocated to long-term wholesale market contracts ($\term_{mcts}$) adheres to the terms of the contract.
Supply-demand balance constraints~\eqref{eq:supply_demand_balance} ensure that the total production equals the total quantities supplied to the markets.
Constraints~\eqref{eq:supply_transportation_balance} ensure that total supply equals the total amount transported to all market locations.
Constraints~\eqref{eq:supply limits} govern supply limits by ensuring that total production is within lower and upper bounds.
Side constraints~\eqref{eq:side_constraints} capture potential mixed-integer requirements through the set $\mc{X}$.  
The remaining constraints are variable bounds.

To handle ``separable'' cross-market price elasticities, i.e., the situation when the volume supplied to market $m_1$ impacts/depresses the price in market $m_2 (\neq m_1)$ as well, one could replace the term $\sum_{k \in \mc{K}_m} P_{mkts} \spot_{mkts}$ in equation~\eqref{eq:profit_in_scenario_s} with a more complex multivariate ``staircase'' representation
\begin{equation}
\sum_{k \in \mc{K}_m} P_{mkts} \spot_{mkts} - \sum_{m' \neq m} \sum_{k' \in \mc{K}_{m'}} \delta_{mm'k'ts} \spot_{m'k'ts}.
\end{equation}
Here, $\delta_{mm'kts}$ denotes the per-unit price reduction in market $m$ due to volumes sold in market $m'$. More sophisticated piecewise linear representations could be used to capture and linearize other nonlinear cross-market price relationships. 
For ease of exposition and due to the additional complexities associated with forecasting these nonlinear relationships, we will henceforth omit this interesting generalization.

\subsection{Conditional value-at-risk formulation}

Transitioning away from a purely risk-neutral mindset, a risk-averse decision maker could consider a two-stage risk-averse stochastic program that maximizes a convex combination of the expected profit and the expected ``tail'' profit or conditional value-at-risk ($\Cvar$) profit \citep{rockafellar:2007}.
Conditional value-at-risk has emerged as an extremely popular risk metric due to its coherency and tractability. As such, it has been widely used in applications of optimization under uncertainty.

To arrive at a risk-averse formulation, let $\lambda \in [0,1]$ be a user-defined parameter governing the weight given to the expected value. When $\lambda=1$, the model reverts to the risk-neutral Formulation~\eqref{model:risk_neutral_basic}, while $\lambda=0$ implies that the decision maker is only concerned with the expected ``tail'' profit.  Let $\alpha \in (0,1)$ denote the risk-aversion parameter (or the $\alpha$-quantile) in the $\Cvar$ calculation. In our maximization setting, $\Cvar_{\alpha}(X)$ denotes the expectation of a random variable (profit) $X$ in the conditional distribution of its $(1-\alpha)$-lower tail, e.g., the average of the lowest $1-\alpha=5\%$ profits.  
These assumptions lead to the following risk-averse MILP formulation (see, e.g., \cite{noyan:2012}):
\begin{subequations} \label{model:cvar}
\begin{align}
\max_{\substack{\v{x}^{\min},v^{\textrm{VaR}},\\\vell,\v{u},\v{x},\v{y},\v{z}}}~& \lambda \sum_{s \in \mc{S}} \pi_{s} z_s + (1-\lambda)\Big[ v^{\textrm{VaR}} - \tfrac{1}{1-\alpha}\sum_{s \in \mc{S}}\pi_s \ell_s \Big]   \\ 
\st~~
	& \ell_{s} \geq v^{\textrm{VaR}} - z_s, \ell_{s} \geq 0 \qquad \forall s \in \mc{S} \\
	& v^{\textrm{VaR}} \in \Re  \\
	& \eqref{eq:profit_in_scenario_s}-\eqref{eq:profit_in_scenario_s_variable_bounds}
\end{align}
\end{subequations}
The additional decision variables $z_s$ and $\ell_{s}$ denote the profit and the nonnegative tail loss in scenario $s$, while $v^{\textrm{VaR}}$ captures the value-at-risk at the $\alpha$ confidence level.

\cmt{There is no clear rule on how to set the $\Cvar$ confidence level $\alpha$ and risk trade-off parameter $\lambda$.  Ultimately, a user must set these parameters just as they would in certain multi-objective optimization settings with weights for various objectives.  In our work, we fix the $\lambda$ parameter and conduct sensitivities over $\alpha$.  The rationale for this choice is that our users had a strong grasp of the risk-neutral approach and were therefore comfortable fixing $\lambda$ while investigating the expected tail profit as a function of $\alpha$.  One could fix $\alpha$ and allow $\lambda$ to vary, but this requires the user to have comfort with the expected tail profit for the selected $\alpha$ level.  When this familiarity has not been established, it is our experience that investigating an array of $\alpha$ values is more beneficial. Later, one could consider varying both $\alpha$ and $\lambda$ in tandem.  In our numerical experiments, we only vary $\alpha$.}

\subsection{Distributionally robust model over a Wasserstein ball} \label{sec:DRO_section}

Although Formulation~\eqref{model:cvar} allows for some degree of risk-aversion, it is not ``ambiguity-averse.'' That is, Formulation~\eqref{model:cvar} assumes that an empirical probability distribution is known with certainty and then attempts to maximize the expected ``tail'' profit (and possibly other terms) with respect to this known distribution.  In contrast, a distributionally robust model assumes that the distribution itself is unknown and belongs to a known family of distributions.  After ``centering'' the family around the empirical data, a DRO approach over a Wasserstein ball maximizes the worst-case expected profit over this chosen family.

To arrive at a tractable DRO formulation over a Wasserstein ball, we describe the price vector $\v{p}$ as a linear function $\v{p}=\v{Q}\vxi + \v{q}$, where $\v{Q}$ and $\v{q}$ are a predefined matrix and vector that map an underlying basis vector $\vxi$ of uncertain parameters to $\v{p}$.  Furthermore, we must specify the Wasserstein radius $\epsilon \in \Re_+$, which defines the permissible deviation from the empirical data. 
Following \citep[Prop. 3]{xie2020tractable}, one can formulate a two-stage DRO model as follows:
\begin{subequations} \label{model:dro_wasserstein}
\begin{alignat}{4}
\max_{\substack{\v{x}^{\min},\\\v{u},\v{x},\v{y},\v{z}}}~& \sum_{s \in \mc{S}} \pi_{s} z_s - \epsilon \sum_{s \in \mc{S}} \pi_{s} ||\v{Q}^{\top}\v{y}_s^{\textrm{Spot}}||_{p^*} \label{objfnc:dro} \\ 
\st~~
	& \eqref{eq:profit_in_scenario_s}-\eqref{eq:profit_in_scenario_s_variable_bounds} 
\end{alignat}
\end{subequations}

Analogous to the $\Cvar$ term $(1-\lambda)\Big[ v^{\textrm{VaR}} - \tfrac{1}{1-\alpha}\sum_{s \in \mc{S}}\pi_s \ell_s \Big]$ in Formulation~\eqref{model:cvar}, the term $\epsilon \sum_{s \in \mc{S}} \pi_{s} ||\v{Q}^{\top}\v{y}_s||_{p^*}$ in \eqref{objfnc:dro} acts as a penalty on the spot allocation.  The larger the $\epsilon$ radius, the larger the penalty on downside volatility. Setting $p=1$ or $p=\infty$, the dual norm $p^*=\infty$ or 1, respectively, and can therefore be represented using linear constraints. We set $p=\infty$ in our computational experiments.

Formally, Formulation~\eqref{model:dro_wasserstein} is 
known as a two-stage distributionally robust linear or mixed-integer linear program under the type-$\infty$ Wasserstein ball.  Both assume a so-called $\infty$-Wasserstein ambiguity set $\mc{P}$. It is known that the $\infty$-Wasserstein ambiguity set offers greater tractability than other Wasserstein ambiguity sets, while still exhibiting attractive convergent properties \citep{xie2020tractable}.

\section{Numerical Results} \label{sec:numerical_results}

This section documents numerical results for two case studies in the PJM electricity market. PJM is a regional transmission organization that coordinates wholesale electricity movement in the mid-Atlantic region of the US. While the first case study assumes a single supplier and single market location, the second case study considers ten locations (nodes) within the PJM market.  There are no mixed-integer constraints in $\mc{X}$, i.e., all of our instances are LPs.  All models were implemented in AIMMS version 4.91.3.6 and solved with Gurobi 10.0. We do not report computational times as all instances could be solved in under a minute. 

\begin{figure}[h!] 
\centering
\includegraphics[width=16cm]{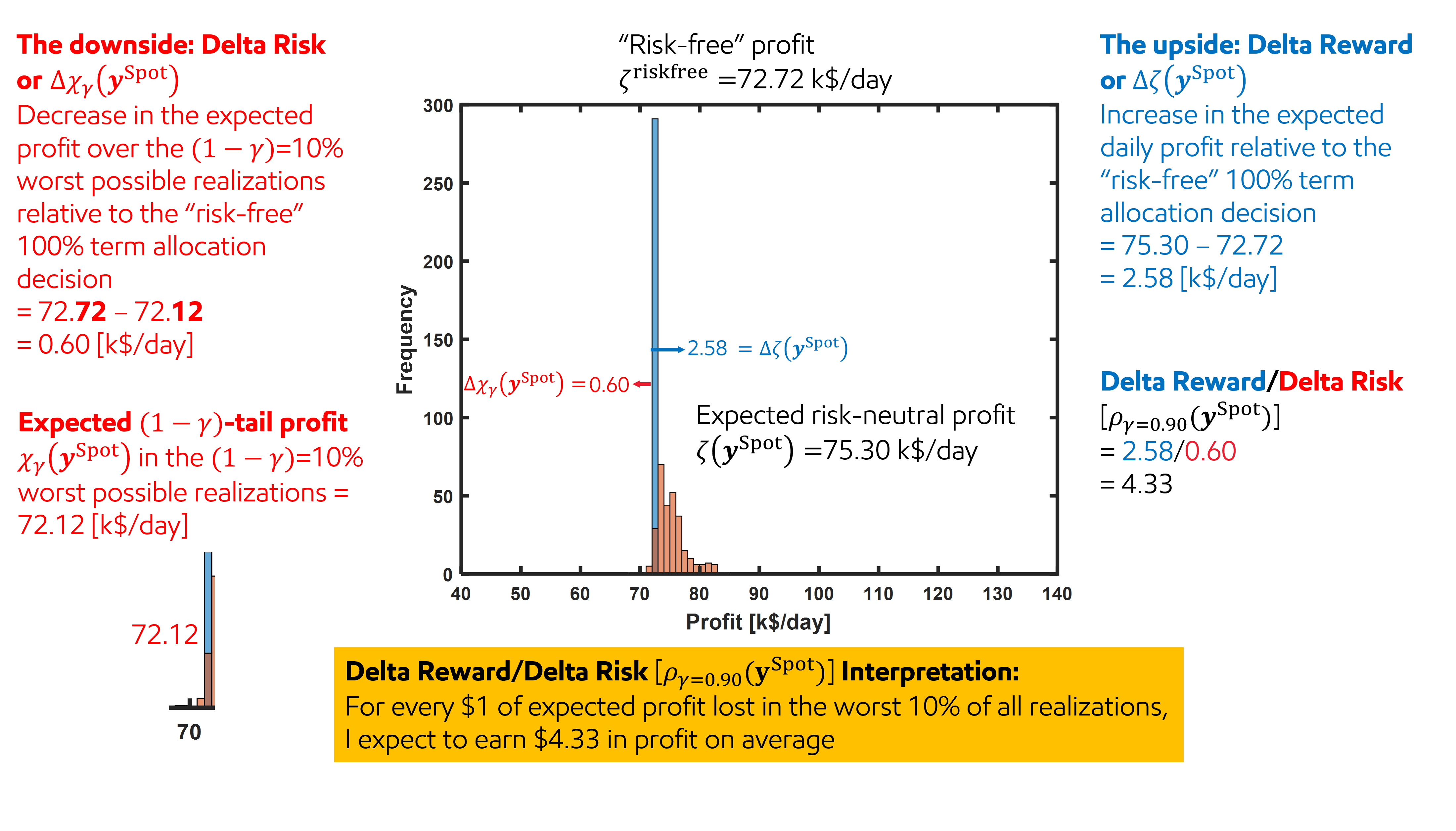}
\caption{Numerical example and explanation of the ``reward-to-risk'' metric $\rho_{\gamma}(\v{y}^{\textrm{Spot}})$ used to compare results. Here $\gamma=0.90$ and $\v{y}^{\textrm{Spot}}$ is a given non-zero spot allocation vector.}
\label{fig:Delta_Reward_over_Delta_Risk_Slide}
\end{figure} 

To compare the three formulations put forth in the previous section,
we need to define several metrics \cmt{(see also Figure~\ref{fig:Delta_Reward_over_Delta_Risk_Slide})}:
\begin{itemize}
\item $\zeta(\v{y}^{\textrm{Spot}})=$ Expected risk-neutral profit given a spot allocation decision $\v{y}^{\textrm{Spot}}$
\item $\chi_{\gamma}(\v{y}^{\textrm{Spot}})=$ Expected $(1-\gamma)$-tail profit ($\Cvar_{\gamma}$) given a spot allocation decision $\v{y}^{\textrm{Spot}}$ and $\gamma \in (0,1)$
\item $\zeta^{\textrm{riskfree}}=$ Expected risk-free profit. Note that $\zeta^{\textrm{riskfree}}=\zeta(\v{0})=\chi_{\gamma}(\v{0}) \quad \forall \gamma \in (0,1)$
\item $\Delta\zeta(\v{y}^{\textrm{Spot}}) =\zeta(\v{y}^{\textrm{Spot}})-\zeta^{\textrm{riskfree}}=$ Expected increase in the risk-neutral profit relative to the risk-free profit given a spot allocation decision $\v{y}^{\textrm{Spot}}$
\item $\Delta\chi_{\gamma}(\v{y}^{\textrm{Spot}}) =|\chi_{\gamma}(\v{y}^{\textrm{Spot}})-\zeta^{\textrm{riskfree}}|=$ Absolute value of the expected decrease in the $(1-\gamma)$-tail profit relative to the risk-free profit given a spot allocation decision $\v{y}^{\textrm{Spot}}$
\item $\rho_{\gamma}(\v{y}^{\textrm{Spot}}) = \Delta\zeta(\v{y}^{\textrm{Spot}}) / \Delta\chi_{\gamma}(\v{y}^{\textrm{Spot}})=$ Change in expected profit per change in risk (expected $(1-\gamma)$-tail profit) as a function of the spot allocation $\v{y}^{\textrm{Spot}}$
\end{itemize}
As a reminder, when we refer to the $(1-\gamma)$-tail profit, we mean the expected value of the worst (lowest) $(1-\gamma) \times 100\%$ profits.  This calculation is akin to sorting the profit per scenario $z_s$ in increasing order, identifying the lowest $(1-\gamma) \times 100\%$, and taking the empirical mean over these values.  Thus, when computing the expectation, we use the finite set of scenarios along with their respective probability.

A word of caution to avoid potential confusion: Throughout this section, we distinguish between two confidence values $\alpha$ and $\gamma$.  The scalar $\alpha$ represents the confidence level used to define $\Cvar_{\alpha}$ in Formulation~\eqref{model:cvar}, which we vary from 0 to 1.  Meanwhile, the scalar $\gamma$ is used in our risk calculation of our reward-to-risk metric $\rho_{\gamma}(\v{y}^{\textrm{Spot}})$. We use two values of $\gamma$, i.e., $\gamma\in\{0.90,0.95\}$ to quantify risk.
Why do we need both $\alpha$ and $\gamma$?
We assume that a decision-maker would like to explore the impact of $\alpha$ and $\epsilon$ on the long-term vs. spot allocation decisions.  Meanwhile, $\gamma$ serves to define a risk metric over which all results can be easily compared. 

\cmt{For both case studies, we apply $k$-means clustering on the raw data and then select the empirical sample closest to the centroid of each cluster. More concretely, we perform a sweep over a range of $k$ values before settling on the ``optimal'' number of clusters using the knee point on a trade-off curve described in \cite{kumaran2021active}. Although $k$-means clustering is prone to favor non-extreme samples, since our main goal is to compare formulations, we opted for this straightforward approach and leave more advanced scenario selection procedures for future work. We provide the scenarios and their associated probabilities for each case study in the Supplementary Information.}

\cmt{
For each of the case studies, we compute the matrix $\v{Q}$ and vector $\v{q}$ for the DRO formulation as follows.
Let $\v{p}_m$ denote the column vector of historical real-time hourly LMPs at market location (node) $m \in \mc{M}$, $\v{p}^{\text{PJM}}$ the vector of observed system-wide LMPs for all of PJM, and $\hat{\v{p}}_m = \v{p}_m-\v{p}^{\text{PJM}}$ the deviation from the system-wide real-time price. 
Let $\hat{\v{P}} = [\hat{\v{p}}_1 \dots \hat{\v{p}}_{|\mc{M}|}]$ be the matrix of LMP deviations from the system-wide price. 
After computing the covariance matrix $\vSigma = \Cov(\hat{\v{P}})$ of LMP deviations, we set $\v{R}=\texttt{chol}(\vSigma)$ where $\v{R}$ is an upper triangular matrix returned by a Cholesky factorization of $\vSigma$ such that $\vSigma = \v{R}^{\top}\v{R}$.  
Note that, for simplicity, $\vSigma \in \Re^{|\mc{M}| \times |\mc{M}|}$ only captures the covariance between markets, not in the time dimension.
We then set $\v{Q}=\v{R}^{\top}$ and $\v{q}=\bar{p}^{\text{PJM}}\v{1}$, where $\bar{p}^{\text{PJM}}$ is the mean of $\v{p}^{\text{PJM}}$ and $\v{1}$ is a column vector of ones.
Finally, we assume that the underlying random basis vector $\vxi \in \Re^{|\mc{M}|}$ discussed in Section~\ref{sec:DRO_section} and \cite{xie2020tractable} satisfies $\E[\vxi]=\v{0}$ and $\Cov(\vxi)=\v{I}$. Then, with these assumptions and letting $\v{p} = \v{Q}\vxi + \v{q}$, we have $\E[\v{p}]=\v{q}$ (or $\E[p_m]=\bar{p}^{\text{PJM}}$ for all $m \in \mc{M}$) and $\Cov(\v{p})=\Cov(\v{Q}\vxi + \v{q})=\Cov(\v{R}^{\top}\vxi + \v{q})=\v{R}^{\top}\v{R}=\vSigma$.
In this sense, we assume ellipsoidal uncertainty in our DRO model.  \cmt{We provide} the data for each case study, including the historical prices and the $\v{Q}$ matrices (recall that the $\v{q}$ vectors do not appear in the DRO formulation), in the Supplementary Information.
}

\subsection{Case Study 1: Single Market Location}
\begin{figure}[h!] 
\centering
\includegraphics[width=9cm]{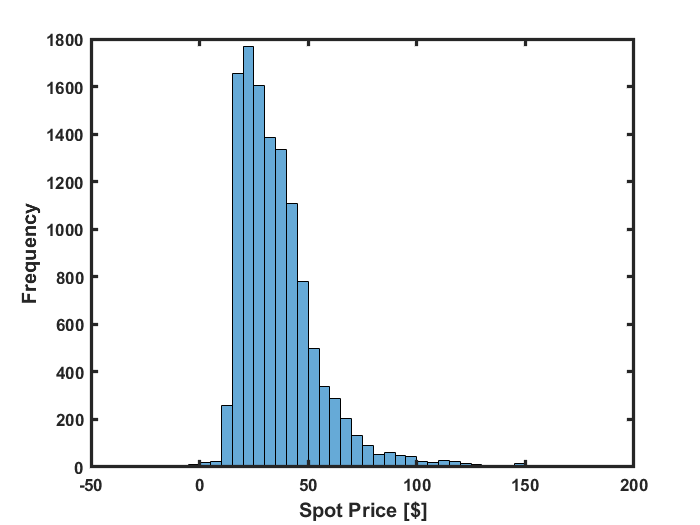}
\caption{Histogram of historical real-time hourly location marginal prices (LMP) at PJM node 48612 from 1 Jan 2021 through 1 May 2022.}
\label{fig:lmp_histogram}
\end{figure}

For simplicity, we assume a single market location $|\mc{M}|=1$, a planning horizon of $|\mc{T}|=365$ days, and $|\mc{S}|=100$ scenarios are available. 
Figure~\ref{fig:lmp_histogram} depicts a truncated histogram of historical real-time hourly location marginal prices (LMPs) at PJM pricing node 48612 from 1 Jan 2021 through 1 May 2022 available at \url{http://dataminer2.pjm.com/feed/rt_hrl_lmps}. The histogram is truncated as prices spiked higher than \$200 in several hours.
There is a single supply cost and no transportation cost.

We assume that $X_{mc}^+ = 0$ for all $c$ (recall $|\mc{M}|=1$), the minimum and maximum supply satisfy $L_t = U_t = 500$ MW for all $t$.
For some perspective, the average capacity of a natural gas-fired combined-cycle power block is roughly 500 MW.
The long-term contract demand curve has individual contracts up to 20 MW each, starting at a price of 38 \$/MWh, decreasing by 1 \$/MW with each subsequent contract, i.e., $W_{mc} = 38 - (c-1)$ and $X_{mc}^{\max} = 20$ MW for all $c=1,\dots,20$. The spot price elasticity curve has steps of width 25 MW, 
while the spot price decreases by 0.2 \$/MWh relative to the nominal parameter value in each step. That is,
$P_{kts} = P_{1ts} - (k-1)0.2$ \$/MWh and $Y_{kts}^{\textrm{Spot}} = 25$ MW for $k=1,\dots,|K|=20$. \cmt{See the righthand side of Figure~\ref{fig:Power_portfolio_opt_problem_statement} for a spot price elasticity curve example when $P_{1ts}=50$ \$/MWh for a given $t$ and $s$. Although not shown in Figure~\ref{fig:Power_portfolio_opt_problem_statement}, the long-term contract demand curve has similar decreasing ``staircase'' structure.}
We set $\lambda=0.01$, thus weighting the $\Cvar$ component of the objective function much more heavily than the expected profit.

\begin{figure}[h!] 
\centering
\includegraphics[width=10cm]{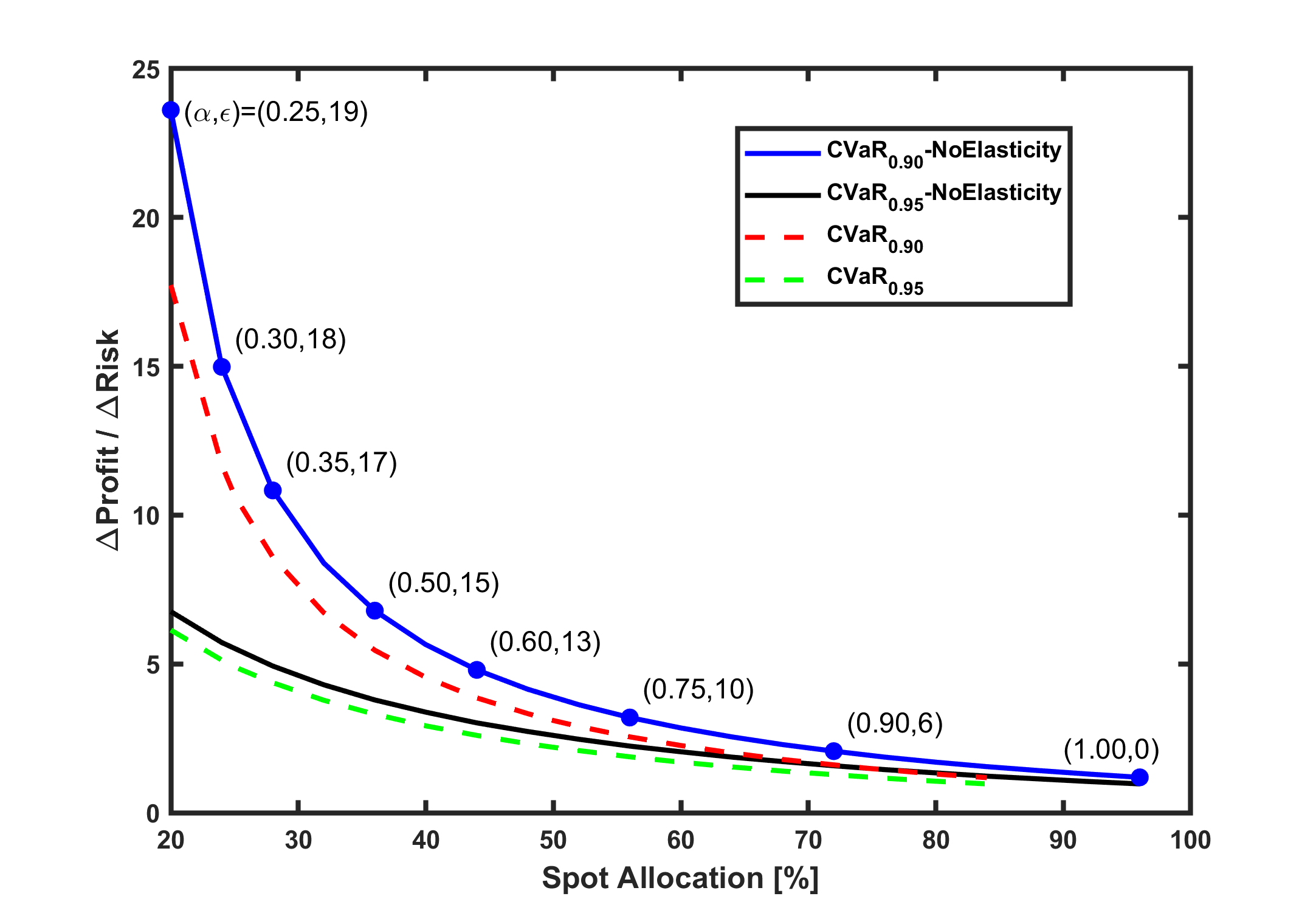}
\caption{Profit vs risk tradeoff curves \cmt{for Case Study 1}. $\Delta$Profit and $\Delta$Risk are \{the expected profit increase\} and \{the absolute value of the $\Cvar_{\gamma}$ decrease\}, respectively, relative to the risk-free allocation of committing all supply to long-term contracts. The $(\alpha,\epsilon)$ values shown next to a subset of points refer to the $\Cvar_{\alpha}$ confidence level and the Wasserstein radius $\epsilon$ that induce the spot allocation on the horizontal axis. Two $\gamma$ values (0.95 and 0.90) are used in the $\Delta$Risk calculation, while $\alpha$ values from 0 to 1 are tested to generate the four curves.}
\label{fig:reward_vs_risk_tradeoff}
\end{figure}

With the metric definitions given above, we can analyze the risk-reward tradeoff associated with allocating more supply to the spot market.
Figure~\ref{fig:reward_vs_risk_tradeoff} depicts the ``change in reward per change in risk'' metric $\rho_{\gamma}(\v{y}^{\textrm{Spot}})$ as a function of the spot allocation.
The risk term in the denominator is measured using both the $\Cvar_{\gamma=0.95}$ and $\Cvar_{\gamma=0.90}$ metrics, which is also reflected in the legend labels. 
Intuitively, as one takes on more risk, the profit per unit risk decreases.
When virtually all supply is allocated to essentially risk-free long-term contracts (see the left side of Figure~\ref{fig:reward_vs_risk_tradeoff}), re-allocating a small amount of supply to the spot market results in a small decrease in the $\Cvar$ and a large increase in expected profit. More colloquially, a practitioner would say that taking on a small amount of risk (measured in terms of $\Cvar$) results in a large expected reward.  Note that the horizontal axis in Figure~\ref{fig:reward_vs_risk_tradeoff} begins with a spot allocation of 20\% because, at a 0\% spot allocation, the ``change in reward per change in risk'' metric $\rho_{\gamma}(\v{y}^{\textrm{Spot}})$ is extremely large due to a small denominator. As more supply is re-allocated to the spot market, however, this profit-per-unit-risk metric tapers off revealing that an additional unit of expected profit is only available by taking on a larger amount of risk.

Some discussion on the usefulness and interpretability of Figure~\ref{fig:reward_vs_risk_tradeoff} is warranted. In our experience, stakeholders do not know the precise value of $\alpha$ for the tail of the $\Cvar$ calculation or the Wasserstein radius $\epsilon$ to use to define their level of risk. 
The tradeoff curves shown in Figure~\ref{fig:reward_vs_risk_tradeoff} guide stakeholders in making actionable decisions by quantifying a range of uncertainty after fixing a value on one of the two axes.  By isolating a spot allocation percentage on the horizontal axis, e.g., a 20\% spot allocation, a decision maker can quantify a range of potential reward-to-risk ratios that result from varying their tolerance for what occurs in the worst possible realizations.  
On the other hand, by isolating a minimum desired reward-to-risk value on the vertical axis, e.g., a ratio of 2-to-1, a decision maker can quantify a range of potential spot market percentages that would yield the desired minimum ratio. In other words, the different risk tolerance curves provide these ranges, which can ultimately give a decision maker more confidence in the actions that they should choose to implement. We have found these tradeoff curves to be more appealing to stakeholders than the traditional risk-reward Pareto curve seen with the classic Markowitz portfolio model.

\cmt{Figure~\ref{fig:reward_vs_risk_tradeoff} also illustrates the impact of including spot price elasticity}.
As stated in the previous section, the presence of spot price elasticity implies that, as the decision maker injects more supply into the market, prices may decrease below their nominal scenario value.
This price (and hence profit) reduction relative to the ``NoElasticity'' curve is most noticeable on the left of Figure~\ref{fig:reward_vs_risk_tradeoff}. This behavior occurs because, when only a small amount of supply is reserved for the spot market, it is able to capitalize on ``high price events,'' which become less profitable as additional supply depresses their prices. Consequently, the numerator $\Delta$Profit is lower relative to the ``NoElasticity'' assumption and the ``change in reward per change in risk'' metric decreases. 

We now contrast the formulations. To this end, it is important to note that different values of the confidence value $\alpha$ and Wasserstein radius $\epsilon$ induce different long-term vs. spot market allocation decisions.
\cmt{Figure~\ref{fig:reward_vs_risk_tradeoff} shows these $(\alpha,\epsilon)$ values above a select number of points}.
As shown by the rightmost point in Figure~\ref{fig:reward_vs_risk_tradeoff}, the risk-neutral Formulation~\eqref{model:risk_neutral_basic} is equivalent to setting the confidence level $\alpha=1$ (recall that, for a random variable $X$, $\Cvar_{1}[X] = \E[X]$ by definition) or the Wasserstein radius $\epsilon=0$. Furthermore, in this example, the risk-neutral approach allocates 96\% of supply to the spot market, but commits 4\% of supply to long-term contracts.  The reason why 4\% is retained is because the expected spot price without elasticity is roughly 37.42 \$/MWh, while the first long-term contract is available for 38 \$/MWh.  
In agreement with intuition, as $\alpha$ decreases and $\epsilon$ increases (i.e., as we move from right to left in Figure~\ref{fig:reward_vs_risk_tradeoff}), risk aversion increases and less supply is allocated to the volatile spot market in favor of greater price stability via long-term contracts. 

\cmt{We generated Figure~\ref{fig:reward_vs_risk_tradeoff} by independently solving Model~\eqref{model:cvar}, while varying the confidence level $\alpha$, and Model~\eqref{model:dro_wasserstein}, while varying the Wasserstein radius $\epsilon$. We then selected $(\alpha,\epsilon)$ pairs that lead to the same spot allocation percentage shown on the horizontal axis.
In general, it is difficult to predict a priori $(\alpha,\epsilon)$ pairs that induce a particular spot allocation percentage. Hence, we solve the two models independently and then generate the curves in a post-processing step.}

Finally, a word on how to interpret the Wasserstein radius $\epsilon$ is in order.
While the interpretation of the confidence level $\alpha$ is well understood as the lower $(1-\alpha)$-tail expectation or largest expected profit in the $(1-\alpha)100\%$ worst scenarios (when maximizing) in the $\Cvar_{\alpha}$ calculation (see, e.g., \cite{rockafellar:2007}), the Wasserstein radius is perhaps less well known.
\citet{ramdas2017wasserstein} provide some intuition. For a one-dimensional random variable, the $\infty$-Wasserstein distance, which is used in our example, between two bounded probability distributions can be interpreted as the maximum distance between the quantile functions of the two distributions. Thus, when $\epsilon=1$, one can interpret the DRO Formulation~\eqref{model:dro_wasserstein} as imposing a limit of $\epsilon=1$ on the maximum distance between the quantile function of the empirical spot price distribution and that of the worst-case spot price distribution.

\subsection{Case Study 2: 10 PJM Nodes}

In our second case study, we consider ten nodes within the PJM electricity market.
Thus, while there is one market, there are ten market locations (i.e., nodes). Table~\ref{tbl:LMP_statistics_CaseStudy2} provides specific information about these nodes including the minimum, median, maximum, and mean hourly LMP used in our data set. Note that the median LMP values are significantly less than the mean values because certain high-price events inflate the mean. 
We list each node's corresponding transmission zone to emphasize that the nodes are geographically dispersed across New Jersey and eastern Pennsylvania.
We assume a planning horizon of $|\mc{T}|=12$ months. \cmt{As described above, we} use $|\mc{S}|=1094$ scenarios after performing \cmt{ $k$-means clustering (with $k=1094$) on the 8760 historical 10-dimensional vectors and then selecting the empirical sample closest to the centroid of each cluster. More concretely, we performed a sweep over a range of $k$ values before settling on $k=1094$ as the preferred number of clusters using the knee point on a trade-off curve explained in \cite{kumaran2021active}.}

\begin{table}[h]									
\centering	
\caption{A list of the 10 PJM nodes used in the second case study along with statistics for their 2022 real-time LMP values. The last row lists these statistics for the system-wide PJM real-time LMP values.}	
\begin{tabular}{ccrrrrrr}									
\toprule
	&	&	& \multicolumn{5}{c}{\textbf{2022 LMP Statistics [\$/MWh]}}	\\	
													\cmidrule(lr){4-8}
Node Name	&	Transmission Zone	&	Node ID	&	Min	&	Median	&	Max	&	Stdev	&	Mean	\\
\midrule
ATLANTIC	&	JCPL	&	2041987701	&	-20.19	&	51.84	&	3584.33	&	110.73	&	71.93	\\
BRIDGEWA	&	PSEG	&	106856853	&	-829.68	&	51.48	&	3530.43	&	111.98	&	68.85	\\
CARDIFF		&	AECO	&	31020665	&	-65.73	&	50.81	&	3580.05	&	109.75	&	70.20	\\
DALESVIL	&	PECO	&	45565905	&	-178.93	&	47.37	&	3462.78	&	111.71	&	65.16	\\
HCS			&	PSEG	&	1097732304	&	-38.82	&	55.70	&	3477.20	&	111.99	&	78.14	\\
HOBOKEN		&	PSEG	&	119118165	&	-309.79	&	53.87	&	3482.89	&	110.59	&	75.20	\\
HUNTERST	&	METED	&	1183231892	&	-19.96	&	56.40	&	3936.18	&	121.69	&	77.49	\\
KINGSLAN	&	PSEG	&	2156109752	&	-106.83	&	53.65	&	3484.85	&	109.75	&	74.39	\\
SAYRECON	&	JCPL	&	119118305	&	-20.17	&	52.27	&	3594.79	&	112.66	&	72.83	\\
WWHARTON	&	JCPL	&	1268571440	&	-70.77	&	53.02	&	3579.92	&	110.60	&	72.55	\\
\midrule
PJM	&		&		&	-20.17	&	57.055	&	3700.00	&	115.58	&	115.58	\\
\bottomrule	
\end{tabular}	
\label{tbl:LMP_statistics_CaseStudy2}	
\end{table}	

\begin{figure}[h!] 
\centering
\includegraphics[width=16cm]{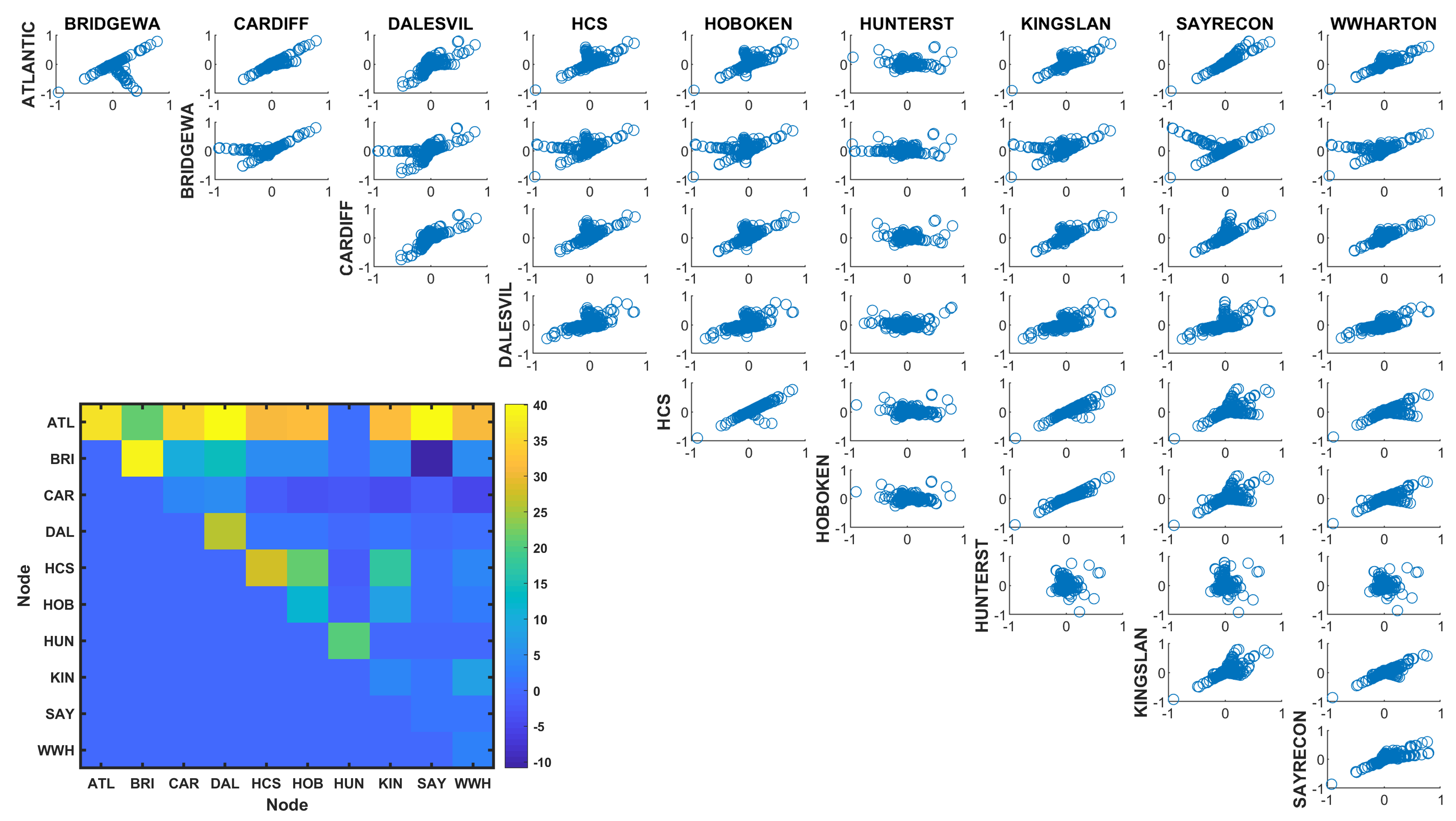}
\caption{Top right: Scatter plots of locational marginal price correlations for 10 PJM market locations (nodes) reveal that prices at certain nodes are strongly positively correlated, while other relationships are more anomalous. Bottom left: Upper triangular $\v{Q}^{\top}$ matrix used in the objective function~\eqref{objfnc:dro} of DRO Model~\eqref{model:dro_wasserstein} for Case Study 2.}
\label{fig:pairplot_10_PJM_Nodes}
\end{figure}

Figure~\ref{fig:pairplot_10_PJM_Nodes} \cmt{(top right)} depicts pairwise correlations between the ten market locations.
\cmt{Specifically, let $\v{p}_n$ denote the vector of 8760 historical real-time hourly LMPs at node $n$, $\v{p}^{\text{PJM}}$ the vector of observed system-wide LMPs for all of PJM, $\hat{\v{p}}_n = \v{p}_n-\v{p}^{\text{PJM}}$ the deviation from the system-wide real-time price, and $\dot{\v{p}}_n$ a normalization of $\hat{\v{p}}_n$ to be within $[-\v{1},\v{1}]$. Then the top right of Figure~\ref{fig:pairplot_10_PJM_Nodes} shows a scatter plot of $(\dot{\v{p}}_{n_1},\dot{\v{p}}_{n_2})$ for each pair of nodes $n_1$ and $n_2$ ($n_2 \neq n_1$).}
While we see several market location pairs that are strongly positively correlated (e.g., Atlantic-Cardiff and Hoboken-Kingslan), we observe other relationships as well.  For example, Hunterst has a nebulous correlation with Kingslan, Sayrecon, and Wwharton. Meanwhile, Bridgewa reveals that there were hours when its 2022 LMP remained flat, while other nodes experienced significant negative prices. 

Similar to our first case study, we assume that $X_{mc}^+ = 0$ for all $m$ and $c$ and that the minimum and maximum supply satisfy $L_t = U_t = 500$ MW for all $t$.
The long-term contract demand curve has individual contracts up to 150 MW each, starting at a price of 62 \$/MWh, decreasing by 1 \$/MW with each subsequent contract, i.e., $W_{c} = 38 - (c-1)$ and $X_{c}^{\max} = 150$ MW for all $c=1,\dots,40$.  We have suppressed the index $m$ for $W_{c}$ and $X_{c}^{\max}$ because we assume that long-term contracts are available to (and can be served from generators in) the entire PJM market.
The spot price elasticity curve has steps of width 50 MW, 
while the spot price decreases by 1 \$/MWh relative to the nominal parameter value in each step. That is,
$P_{kts} = P_{1ts} - (k-1)$ \$/MWh and $Y_{kts}^{\textrm{Spot}} = 50$ MW for $k=1,\dots,|K|=10$. \cmt{Once again, see the righthand side of Figure~\ref{fig:Power_portfolio_opt_problem_statement} for a spot price elasticity curve example when $P_{1ts}=50$ \$/MWh for a given $t$ and $s$.}
We set $\lambda=0.10$, thus favoring the $\Cvar$ contribution of the objective function more than the expected profit.
There is a single supply cost and no transportation cost.


\cmt{As described above}, given the 8760 prices for each PJM node in the year 2022, we compute the matrix $\v{Q} \in \Re^{10 \times 10}$ needed in the DRO Formulation~\eqref{model:dro_wasserstein} as follows.  We compute the empirical covariance matrix $\vSigma$ corresponding to normalized nodal prices.  Specifically, letting $\v{p}_n$ denote the vector of observed LMPs at node $n$ and $\v{p}^{\text{PJM}}$ denote the vector of observed LMPs for all of PJM, we obtain $\vSigma$ by computing the covariance of the vectors $\v{p}_n-\v{p}^{\text{PJM}}$ for all $n$.  
We then take $\v{Q}$ to be the lower triangular matrix returned via Cholesky decomposition of $\vSigma$. \cmt{(Most linear algebra codes, e.g., Matlab's \texttt{chol} function, return an upper triangular matrix $\v{R}$ such that $\v{R}^{\top}\v{R}=\vSigma$; we take $\v{Q}=\v{R}^{\top}$.)} This amounts to an ellipsoidal uncertainty set for spot prices in our DRO model \citep{matthews2018generalized}.
\cmt{Figure~\ref{fig:pairplot_10_PJM_Nodes} (bottom left) depicts a heat map of the $\v{Q}^{\top}$ matrix appearing in the objective function~\eqref{objfnc:dro} of DRO Model~\eqref{model:dro_wasserstein}.}

\begin{figure}[h!] 
\centering
\includegraphics[width=10cm]{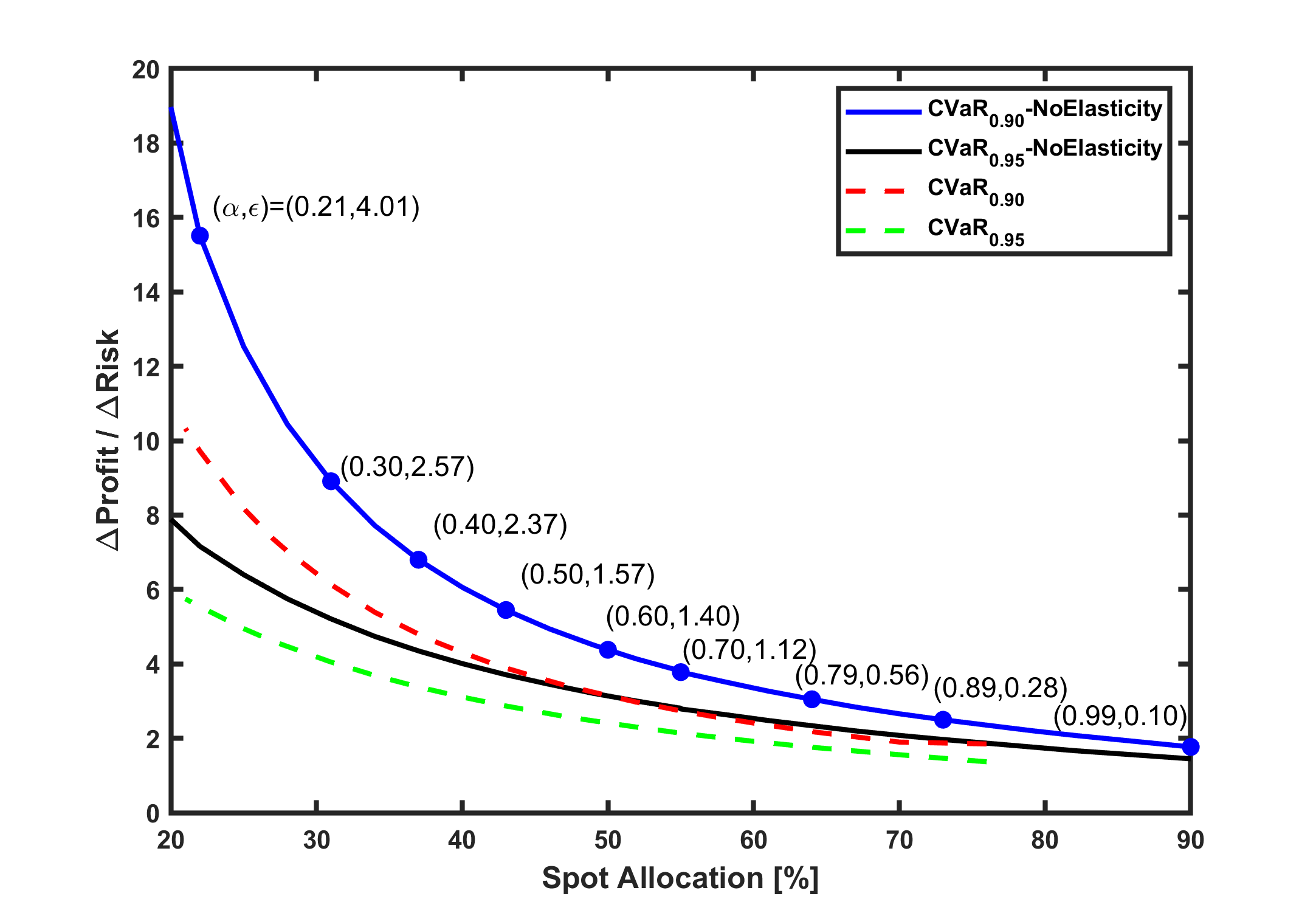}
\caption{Profit vs risk tradeoff curves for Case Study 2.}
\label{fig:reward_vs_risk_tradeoff2}
\end{figure}

Figure~\ref{fig:reward_vs_risk_tradeoff2} depicts the ``change in reward per change in risk'' metric $\rho_{\gamma}(\v{y}^{\textrm{Spot}})$ as a function of the spot allocation for our 10-node case study.
As before, \cmt{we measure} the risk term in the denominator using both the $\Cvar_{\gamma=0.95}$ and $\Cvar_{\gamma=0.90}$ metrics. 
The figure closely mirrors the behavior shown in Figure~\ref{fig:reward_vs_risk_tradeoff} where the unit profit per unit risk decreases as one takes on more risk. One difference between the two is that, when \cmt{we assume} no spot price elasticity, the risk-neutral solution yields a roughly 90\% allocation to the spot market in the present case study, whereas it was closer to 96\% in Case Study 1.

\begin{figure*}[t!]
    \centering
    \begin{subfigure}[t]{0.5\textwidth}
        \centering
        \includegraphics[height=6.5cm]{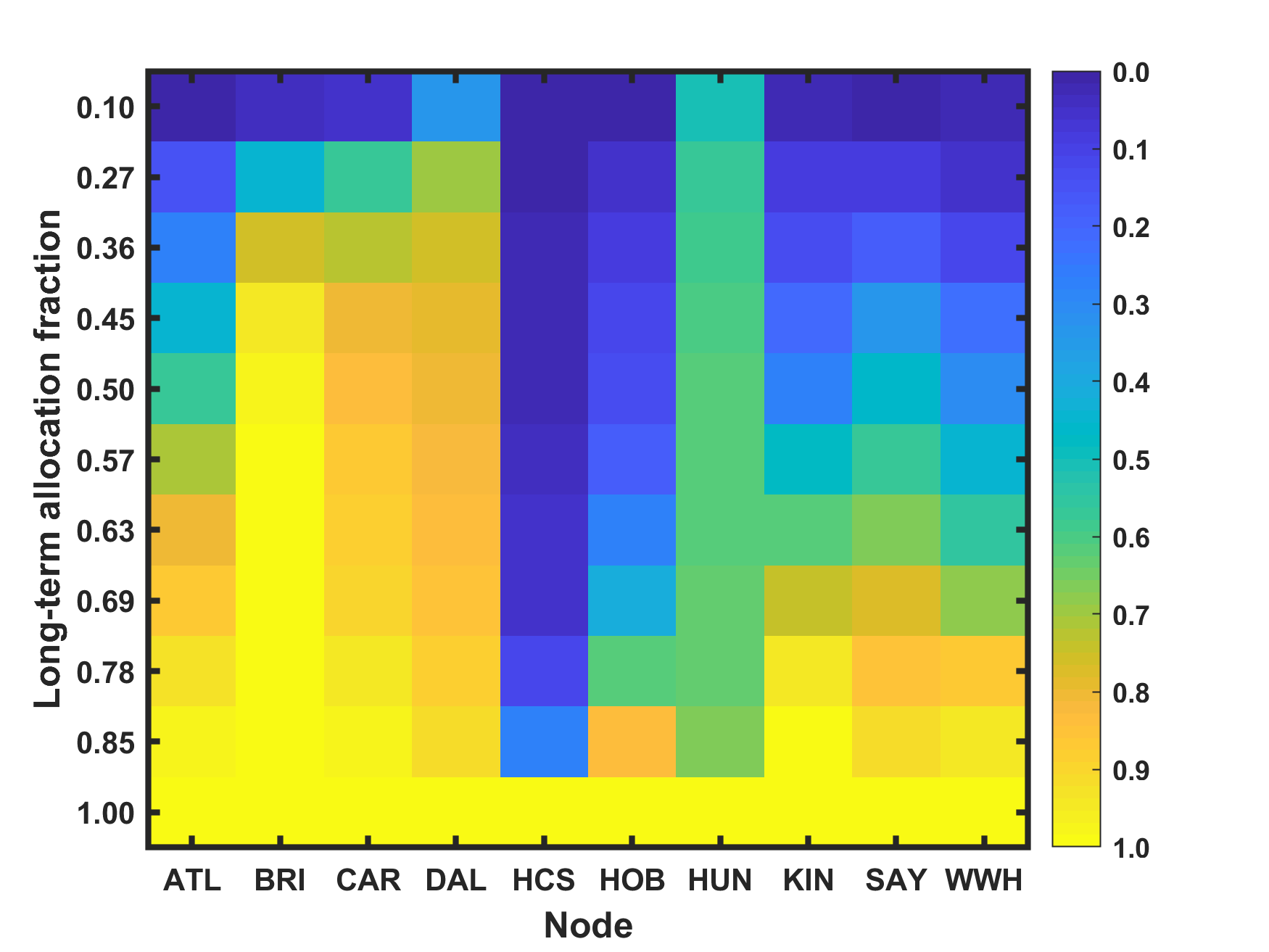}
        \caption{$\Cvar$ Model~\eqref{model:cvar}}
    \end{subfigure}%
    ~ 
    \begin{subfigure}[t]{0.5\textwidth}
        \centering
        \includegraphics[height=6.5cm]{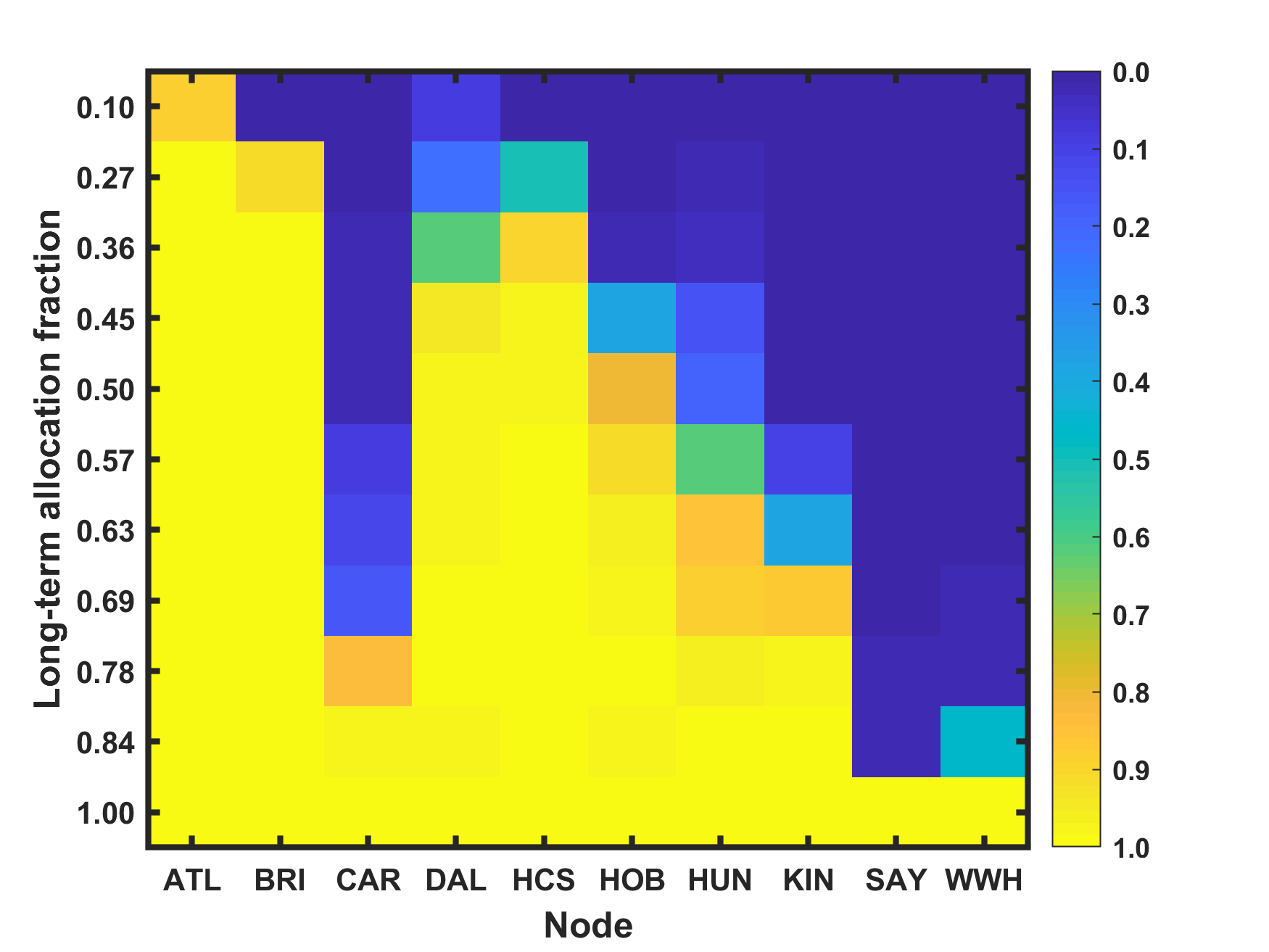}
        \caption{DRO Model~\eqref{model:dro_wasserstein}}
    \end{subfigure}
    \caption{Mean long-term wholesale allocation fraction for Case Study 2 with no spot price elasticity.}
    \label{fig:CaseStudy2_Nodal_Results}
\end{figure*}

\cmt{Figure~\ref{fig:CaseStudy2_Nodal_Results} } shows the mean long-term wholesale allocation fraction for each of the ten nodes as \cmt{we vary} $\alpha$ (used to compute $\Cvar$ in Formulation~\eqref{model:cvar}) and the Wasserstein radius $\epsilon$ (used in Formulation~\eqref{model:dro_wasserstein}). \cmt{We assume} no spot price elasticity. 
\cmt{Since the spot allocation to each node can change by scenario, we show the average long-term wholesale allocation fraction, a value between 0 and 1, to show how certain nodes are more risk-averse (i.e., prefer long-term wholesale contracts) than others.} 
Consistent with the two solid lines (labeled with the suffix ``NoElasticity'') shown in Figure~\ref{fig:reward_vs_risk_tradeoff2}, it is clear that as $1-\alpha$ and $\epsilon$ increase so too does the mean wholesale allocation at each node. 

\cmt{Looking first at the nodal allocation results of the $\Cvar$ Model~\eqref{model:cvar} in Figure~\ref{fig:CaseStudy2_Nodal_Results}(a)}, with the highest mean LMP and second highest median LMP \cmt{(see Table~\ref{tbl:LMP_statistics_CaseStudy2})}, HCS prefers to allocate supply to the spot market more so than any other node.  Meanwhile, with the smallest mean and median LMP, Dalesvil allocates the most supply to the long-term contracts.  Upon first inspection, arguably the most ``non-intuitive'' observation is that Hunterst also prefers to allocate the majority of its supply to long-term contracts when it possesses the largest median and maximum LMP as well as the second largest mean LMP.  However, upon further inspection, the scatter plots in Figure~\ref{fig:pairplot_10_PJM_Nodes} indicate that, when Hunterst's LMPs are large (far to the right), the LMPs at the other nodes remain relatively flat. For many other node pairs, there is a strong positive correlation amongst prices.  As a consequence, despite the large median and mean LMP values, the $\Cvar$ Model~\eqref{model:cvar} prefers to allocate Hunterst's supply to long-term contracts so that it can take advantage of the subdued correlation in LMP values.

\cmt{
A visual inspection of the heat maps in Figure~\ref{fig:CaseStudy2_Nodal_Results} clearly reveals that DRO Model~\eqref{model:dro_wasserstein} yields different results from $\Cvar$ Model~\eqref{model:cvar} at the nodal level. Since the two models share the same structural constraints~\eqref{eq:profit_in_scenario_s}--\eqref{eq:profit_in_scenario_s_variable_bounds}, the reason for the different results lies in the objective function~\eqref{objfnc:dro}.  Specifically,   
the penalty term $\epsilon \sum_{s \in \mc{S}} \pi_{s} ||\v{Q}^{\top}\v{y}_s||_{p^*}$ elucidates the manner in which the estimated $\v{Q}$ matrix interacts with spot allocation decisions.
We refer the reader once again to the heat map of the $\v{Q}^{\top}$ matrix shown in the bottom left of Figure~\ref{fig:pairplot_10_PJM_Nodes}. Bridgewa and Atlantic have the largest diagonal component with $Q_{\text{BRI},\text{BRI}}=38.82$ and $Q_{\text{ATL},\text{ATL}}=36.61$, followed by HCS (27.64), Dalesvil (25.95), Hunterst (20.60), and Hoboken (11.92). In the DRO model, as $\epsilon$ increases (i.e., as the penalty on spot allocations is increased), these nodes witness the least spot allocation percentage almost precisely in the decreasing order of their $Q_{m,m}$ value.  The two exceptions are: Atlantic and Bridgewa swap places as do Hoboken and Hunterst. Atlantic exhibits a higher degree of long-term allocation than Bridgewa in Figure~\ref{fig:CaseStudy2_Nodal_Results}(b) because Atlantic has larger $Q_{\text{ATL},m}$ entries in the $\v{Q}^{\top}$ matrix (see the top row of $\v{Q}^{\top}$) relative to Bridgewa's $Q_{\text{BRI},m}$ entries. A similar observation holds for Hoboken and Hunterst. In summary, we see that the $\v{Q}$ matrix, which can loosely be interpreted as the square root of the covariance matrix $\vSigma$, governs the weight of the spot allocation decisions in DRO Model~\eqref{model:dro_wasserstein} much like the covariance matrix governs the behavior of a traditional Markowitz mean-variance portfolio optimization model \citep{luenberger2014investment}.}

\cmt{The different nodal results ultimately lead to the question: Which model is best, i.e., gives the most trustworthy and useful recommendations? The answer is: It depends.  Both models are governed by fundamental assumptions about the data.  Recall that we fed historical spot price data to both risk-averse models. While $\Cvar$ Model~\eqref{model:cvar} relied on this data to compute the expected tail profit in unfavorable scenarios, DRO Model~\eqref{model:dro_wasserstein} required us to state a $\v{Q}$ matrix vector to explicitly define an affine relationship between the uncertain primitive vector $\vxi$ and the uncertain price vector $\v{p}$.  After estimating $\v{Q}$ as is traditionally done in financial applications (see, e.g., \cite{elton2014modern}), we instantiated and solved the DRO model to determine long-term and spot allocation decisions.  In machine learning vernacular, we solved both models on ``training'' data and should perform an independent experiment on ``test'' data to validate the performance of each.  This validation step is out of scope for this work as our primary goal was to introduce the two risk-averse models for this particular allocation problem.}

\cmt{As a final note, while we empirically found that there is almost always an $(\alpha,\epsilon)$ pair that lead to the same aggregate long-term wholesale allocation, Figure~\ref{fig:CaseStudy2_Nodal_Results} reveals that such a pair does not always exist.  Specifically, one can inspect the vertical axis and note that the second largest long-term allocation fraction is 0.85 in $\Cvar$ Model~\eqref{model:cvar}, but 0.84 in DRO Model~\eqref{model:dro_wasserstein}.  We could not find an $\epsilon$ that induces an aggregate long-term allocation fraction of 0.85.}

\section{Conclusions and future research directions} \label{sec:conclusions}

In this work, we have attempted to make a case for applying DRO to the rather generic problem of balancing supply allocation to long-term contracts where price stability reigns versus spot markets where volatility can lead to large gains and losses. Focusing on a Genco's medium-term planning problem, we presented risk-neutral, risk-averse ($\Cvar$), and ambiguity-averse (DRO) formulations to address it. These formulations also included price elasticity components atypical in the majority of price-taker models.  We then demonstrated how a risk-averse and an ambiguity-averse approach converge to the same decisions depending on the parameter values chosen. As always, the final long-term vs. spot market allocation decision depends on the risk/ambiguity aversion of the decision maker. 

Figures~\ref{fig:reward_vs_risk_tradeoff} and \ref{fig:reward_vs_risk_tradeoff2} reveal that the $\Cvar$ and DRO approaches both produce the same ``change in reward per change in risk'' tradeoff curves, just in a different manner.  Whereas the $\Cvar$ approach requires the user to vary $\alpha$, the probability level governing the conditional expectation that they would like to use to quantify worst-case profits, the DRO approach requires adjustments to the Wasserstein radius $\epsilon$, which governs the magnitude of ambiguity assumed in the data. While the two approaches yield the same \cmt{aggregate} tradeoff curves, it is important to mention that key stakeholders may not have equal understanding in the philosophy behind how the results were generated.  In our experience, it has taken time to explain, even to sophisticated users, how to interpret $\alpha$ and $\epsilon$ in the context of their problem and their data. Ultimately, one interpretation may give users more confidence, which is critical to having them use the decision support models for guidance.

As for future research directions, it would be interesting to consider simultaneous uncertainty in the objective function and the constraints. \cmt{This structure would arise if one were to model transmission limits on the amount of supply that can be allocated to a particular market in certain scenarios.}
Additionally, one could explore price setter behavior or game-theoretic models in which suppliers could exert market power and could therefore act as Nash-Cournot players.
Binary decisions could, of course, be included to model fixed costs associated with production and/or transportation.  Such additions would lead to risk-averse mixed-integer linear optimization problems. 
\cmt{One could explore more advanced scenario selection methods to determine the scenarios to feed each of the risk-averse models. After exploring different techniques for generating the $\v{Q}$ matrix used in DRO Model~\eqref{model:dro_wasserstein}, it would also be beneficial to perform a comprehensive validation experiment to determine which risk-averse model is the most reliable.}

Finally, although we have focused our investigation on a simplified electricity sector application, we are confident that there is potential to apply a DRO framework to other process systems engineering applications. \cmt{The most obvious analogs to the applications mentioned at the end of Section~\ref{sec:electricity_applications} are those involving multi-scale integration. For example, problems that include sizing decisions (e.g., of separation/production units, inventory tanks, energy storage units, and other critical infrastructure) and detailed operational/dispatch decisions could benefit from data-driven approaches under uncertainty. For these problems, sizing decisions are roughly tantamount to the long-term contract decisions studied in this work. Demand response contracts for energy intensive processes would be another candidate.}

\section*{Acknowledgments}
We wish to thank Vignesh Subramanian and several other colleagues for their insightful suggestions that helped improve the quality of this work. \cmt{We also thank two anonymous referees whose astute observations improved the quality of this manuscript.}

\bibliographystyle{plainnat}
\bibliography{wholesale_vs_spot_refs,cvar_refs}%

\begin{thebibliography}{50}
\providecommand{\natexlab}[1]{#1}
\providecommand{\url}[1]{\texttt{#1}}
\expandafter\ifx\csname urlstyle\endcsname\relax
  \providecommand{\doi}[1]{doi: #1}\else
  \providecommand{\doi}{doi: \begingroup \urlstyle{rm}\Url}\fi

\bibitem[Bhattacharyya et~al.(2011)Bhattacharyya, Turton, and
  Zitney]{bhattacharyya2011steady}
Debangsu Bhattacharyya, Richard Turton, and Stephen~E Zitney.
\newblock Steady-state simulation and optimization of an integrated
  gasification combined cycle power plant with {CO}2 capture.
\newblock \emph{Industrial \& Engineering Chemistry Research}, 50\penalty0
  (3):\penalty0 1674--1690, 2011.

\bibitem[Bienstock et~al.(2020)Bienstock, Escobar, Gentile, and
  Liberti]{bienstock2020mathematical}
Dan Bienstock, Mauro Escobar, Claudio Gentile, and Leo Liberti.
\newblock Mathematical programming formulations for the alternating current
  optimal power flow problem.
\newblock \emph{4OR}, 18\penalty0 (3):\penalty0 249--292, 2020.

\bibitem[Cao et~al.(2016)Cao, Swartz, and Flores-Cerrillo]{cao2016optimal}
Yanan Cao, Christopher~LE Swartz, and Jesus Flores-Cerrillo.
\newblock Optimal dynamic operation of a high-purity air separation plant under
  varying market conditions.
\newblock \emph{Industrial \& engineering chemistry research}, 55\penalty0
  (37):\penalty0 9956--9970, 2016.

\bibitem[Chen and Paschalidis(2018)]{chen2018robust}
Ruidi Chen and Ioannis~C Paschalidis.
\newblock A robust learning approach for regression models based on
  distributionally robust optimization.
\newblock \emph{Journal of Machine Learning Research}, 19\penalty0 (13), 2018.

\bibitem[Conejo et~al.(2008)Conejo, Garcia-Bertrand, Carrion, Caballero, and
  de~Andres]{conejo2008optimal}
Antonio~J Conejo, Raquel Garcia-Bertrand, Miguel Carrion, {\'A}ngel Caballero,
  and Antonio de~Andres.
\newblock Optimal involvement in futures markets of a power producer.
\newblock \emph{IEEE Transactions on Power Systems}, 23\penalty0 (2):\penalty0
  703--711, 2008.

\bibitem[Delage and Ye(2010)]{delage2010distributionally}
Erick Delage and Yinyu Ye.
\newblock Distributionally robust optimization under moment uncertainty with
  application to data-driven problems.
\newblock \emph{Operations research}, 58\penalty0 (3):\penalty0 595--612, 2010.

\bibitem[Dowling and Zavala(2018)]{dowling2018economic}
Alexander~W Dowling and Victor~M Zavala.
\newblock Economic opportunities for industrial systems from frequency
  regulation markets.
\newblock \emph{Computers \& Chemical Engineering}, 114:\penalty0 254--264,
  2018.

\bibitem[Dowling et~al.(2017)Dowling, Kumar, and Zavala]{dowling2017multi}
Alexander~W Dowling, Ranjeet Kumar, and Victor~M Zavala.
\newblock A multi-scale optimization framework for electricity market
  participation.
\newblock \emph{Applied Energy}, 190:\penalty0 147--164, 2017.

\bibitem[Elton et~al.(2014)Elton, Gruber, Brown, and
  Goetzmann]{elton2014modern}
E.J. Elton, M.J. Gruber, S.J. Brown, and W.N. Goetzmann.
\newblock \emph{Modern Portfolio Theory and Investment Analysis}.
\newblock Wiley, 2014.

\bibitem[{ERCOT}(May 2020)]{ERCOT2020}
{ERCOT}.
\newblock Impact of system constraints on pricing in {ERCOT}, May 2020.
\newblock URL
  \url{https://www.ercot.com/files/docs/2020/05/28/Pricing_in_ERCOT_FINAL.pdf}.

\bibitem[Fanzeres et~al.(2014)Fanzeres, Street, and
  Barroso]{fanzeres2014contracting}
Bruno Fanzeres, Alexandre Street, and Luiz~Augusto Barroso.
\newblock Contracting strategies for renewable generators: A hybrid stochastic
  and robust optimization approach.
\newblock \emph{IEEE Transactions on Power Systems}, 30\penalty0 (4):\penalty0
  1825--1837, 2014.

\bibitem[Gabriel et~al.(2012)Gabriel, Conejo, Fuller, Hobbs, and
  Ruiz]{gabriel2012complementarity}
Steven~A Gabriel, Antonio~J Conejo, J~David Fuller, Benjamin~F Hobbs, and
  Carlos Ruiz.
\newblock \emph{Complementarity modeling in energy markets}, volume 180.
\newblock Springer Science \& Business Media, 2012.

\bibitem[Gao et~al.(2019)Gao, Ning, and You]{gao2019data}
Jiyao Gao, Chao Ning, and Fengqi You.
\newblock Data-driven distributionally robust optimization of shale gas supply
  chains under uncertainty.
\newblock \emph{AIChE Journal}, 65\penalty0 (3):\penalty0 947--963, 2019.

\bibitem[Gao and Kleywegt(2016)]{gao2016distributionally}
Rui Gao and Anton~J Kleywegt.
\newblock Distributionally robust stochastic optimization with {W}asserstein
  distance.
\newblock \emph{arXiv preprint arXiv:1604.02199}, 2016.

\bibitem[Guo et~al.(2022)Guo, Bodur, and Papageorgiou]{guo2022generation}
Cheng Guo, Merve Bodur, and Dimitri~J Papageorgiou.
\newblock Generation expansion planning with revenue adequacy constraints.
\newblock \emph{Computers \& Operations Research}, 142:\penalty0 105736, 2022.

\bibitem[Hanasusanto and Kuhn(2018)]{hanasusanto2018conic}
Grani~A Hanasusanto and Daniel Kuhn.
\newblock Conic programming reformulations of two-stage distributionally robust
  linear programs over {W}asserstein balls.
\newblock \emph{Operations Research}, 66\penalty0 (3):\penalty0 849--869, 2018.

\bibitem[Kaye et~al.(1990)Kaye, Outhred, and Bannister]{kaye1990forward}
RJ~Kaye, HR~Outhred, and CH~Bannister.
\newblock Forward contracts for the operation of an electricity industry under
  spot pricing.
\newblock \emph{IEEE Transactions on Power Systems}, 5\penalty0 (1):\penalty0
  46--52, 1990.

\bibitem[Kazempour et~al.(2010)Kazempour, Conejo, and
  Ruiz]{kazempour2010strategic}
S~Jalal Kazempour, Antonio~J Conejo, and Carlos Ruiz.
\newblock Strategic generation investment using a complementarity approach.
\newblock \emph{IEEE transactions on power systems}, 26\penalty0 (2):\penalty0
  940--948, 2010.

\bibitem[Kelley et~al.(2018)Kelley, Baldick, and Baldea]{kelley2018demand}
Morgan~T Kelley, Ross Baldick, and Michael Baldea.
\newblock Demand response operation of electricity-intensive chemical processes
  for reduced greenhouse gas emissions: application to an air separation unit.
\newblock \emph{ACS Sustainable Chemistry \& Engineering}, 7\penalty0
  (2):\penalty0 1909--1922, 2018.

\bibitem[Kirschen and Strbac(2018)]{kirschen2018fundamentals}
Daniel~S Kirschen and Goran Strbac.
\newblock \emph{Fundamentals of power system economics}.
\newblock John Wiley \& Sons, 2018.

\bibitem[Kumaran et~al.(2021)Kumaran, Papageorgiou, Takac, Lueg, and
  Sahinidis]{kumaran2021active}
Krishnan Kumaran, Dimitri~J Papageorgiou, Martin Takac, Laurens Lueg, and
  Nicolas~V Sahinidis.
\newblock Active metric learning for supervised classification.
\newblock \emph{Computers \& Chemical Engineering}, 144:\penalty0 107132, 2021.

\bibitem[Kwon et~al.(2006)Kwon, Scott~Rogers, and Yau]{kwon2006stochastic}
Roy~H Kwon, J~Scott~Rogers, and Sheena Yau.
\newblock Stochastic programming models for replication of electricity forward
  contracts for industry.
\newblock \emph{Naval Research Logistics (NRL)}, 53\penalty0 (7):\penalty0
  713--726, 2006.

\bibitem[Lara et~al.(2018)Lara, Mallapragada, Papageorgiou, Venkatesh, and
  Grossmann]{lara2018deterministic}
Cristiana~L Lara, Dharik~S Mallapragada, Dimitri~J Papageorgiou, Aranya
  Venkatesh, and Ignacio~E Grossmann.
\newblock Deterministic electric power infrastructure planning: Mixed-integer
  programming model and nested decomposition algorithm.
\newblock \emph{European Journal of Operational Research}, 271\penalty0
  (3):\penalty0 1037--1054, 2018.

\bibitem[Liu and Yuan(2021)]{liu2021multistage}
Botong Liu and Zhihong Yuan.
\newblock Multistage distributionally robust design of a renewable source
  processing network under uncertainty.
\newblock \emph{Industrial \& Engineering Chemistry Research}, 60\penalty0
  (21):\penalty0 7883--7903, 2021.

\bibitem[Lorca and Prina(2014)]{lorca2014power}
{\'A}lvaro Lorca and Jos{\'e} Prina.
\newblock Power portfolio optimization considering locational electricity
  prices and risk management.
\newblock \emph{Electric Power Systems Research}, 109:\penalty0 80--89, 2014.

\bibitem[Luenberger(2014)]{luenberger2014investment}
D.G. Luenberger.
\newblock \emph{Investment Science}.
\newblock Oxford University Press, 2014.

\bibitem[Matthews et~al.(2018)Matthews, Guzman, and
  Floudas]{matthews2018generalized}
Logan~R Matthews, Yannis~A Guzman, and Christodoulos~A Floudas.
\newblock Generalized robust counterparts for constraints with bounded and
  unbounded uncertain parameters.
\newblock \emph{Computers \& Chemical Engineering}, 116:\penalty0 451--467,
  2018.

\bibitem[Mayer and Tr{\"u}ck(2018)]{mayer2018electricity}
Klaus Mayer and Stefan Tr{\"u}ck.
\newblock Electricity markets around the world.
\newblock \emph{Journal of Commodity Markets}, 9:\penalty0 77--100, 2018.

\bibitem[Mitra et~al.(2013)Mitra, Sun, and Grossmann]{mitra2013optimal}
Sumit Mitra, Lige Sun, and Ignacio~E Grossmann.
\newblock Optimal scheduling of industrial combined heat and power plants under
  time-sensitive electricity prices.
\newblock \emph{Energy}, 54:\penalty0 194--211, 2013.

\bibitem[Mohajerin~Esfahani and Kuhn(2018)]{esfahani2018dro}
Peyman Mohajerin~Esfahani and Daniel Kuhn.
\newblock Data-driven distributionally robust optimization using the
  {W}asserstein metric: Performance guarantees and tractable reformulations.
\newblock \emph{Mathematical Programming}, 171\penalty0 (1):\penalty0 115--166,
  2018.

\bibitem[Noyan(2012)]{noyan:2012}
Nilay Noyan.
\newblock Risk-averse two-stage stochastic programming with an application to
  disaster management.
\newblock \emph{Computers \& Operations Research}, 39\penalty0 (3):\penalty0
  541--559, 2012.

\bibitem[Palys et~al.(2018)Palys, Allman, and Daoutidis]{palys2018exploring}
Matthew~J Palys, Andrew Allman, and Prodromos Daoutidis.
\newblock Exploring the benefits of modular renewable-powered ammonia
  production: A supply chain optimization study.
\newblock \emph{Industrial \& Engineering Chemistry Research}, 58\penalty0
  (15):\penalty0 5898--5908, 2018.

\bibitem[Pattison et~al.(2016)Pattison, Touretzky, Johansson, Harjunkoski, and
  Baldea]{pattison2016optimal}
Richard~C Pattison, Cara~R Touretzky, Ted Johansson, Iiro Harjunkoski, and
  Michael Baldea.
\newblock Optimal process operations in fast-changing electricity markets:
  framework for scheduling with low-order dynamic models and an air separation
  application.
\newblock \emph{Industrial \& Engineering Chemistry Research}, 55\penalty0
  (16):\penalty0 4562--4584, 2016.

\bibitem[Pineda and Conejo(2012)]{pineda2012managing}
Salvador Pineda and Antonio~J Conejo.
\newblock Managing the financial risks of electricity producers using options.
\newblock \emph{Energy economics}, 34\penalty0 (6):\penalty0 2216--2227, 2012.

\bibitem[Rahimian and Mehrotra(2019)]{rahimian2019droreview}
Hamed Rahimian and Sanjay Mehrotra.
\newblock Distributionally robust optimization: A review.
\newblock \emph{arXiv preprint arXiv:1908.05659}, 2019.

\bibitem[Ralph and Smeers(2006)]{ralph2006EPECs}
Daniel Ralph and Yves Smeers.
\newblock {EPEC}s as models for electricity markets.
\newblock In \emph{2006 IEEE PES Power Systems Conference and Exposition},
  pages 74--80, 2006.

\bibitem[Ramdas et~al.(2017)Ramdas, Trillos, and Cuturi]{ramdas2017wasserstein}
Aaditya Ramdas, Nicol{\'a}s~Garc{\'\i}a Trillos, and Marco Cuturi.
\newblock On {W}asserstein two-sample testing and related families of
  nonparametric tests.
\newblock \emph{Entropy}, 19\penalty0 (2):\penalty0 47, 2017.

\bibitem[Risbeck et~al.(2017)Risbeck, Maravelias, Rawlings, and
  Turney]{risbeck2017mixed}
Michael~J Risbeck, Christos~T Maravelias, James~B Rawlings, and Robert~D
  Turney.
\newblock A mixed-integer linear programming model for real-time cost
  optimization of building heating, ventilation, and air conditioning
  equipment.
\newblock \emph{Energy and Buildings}, 142:\penalty0 220--235, 2017.

\bibitem[Rockafellar(2007)]{rockafellar:2007}
R~Tyrrell Rockafellar.
\newblock Coherent approaches to risk in optimization under uncertainty.
\newblock \emph{Tutorials in Operations Research}, 3:\penalty0 38--61, 2007.

\bibitem[Sen et~al.(2006)Sen, Yu, and Genc]{sen2006stochastic}
Suvrajeet Sen, Lihua Yu, and Talat Genc.
\newblock A stochastic programming approach to power portfolio optimization.
\newblock \emph{Operations Research}, 54\penalty0 (1):\penalty0 55--72, 2006.

\bibitem[Shahidehpour and Alomoush(2017)]{shahidehpour2017restructured}
Mohammad Shahidehpour and Muwaffaq Alomoush.
\newblock \emph{Restructured electrical power systems: Operation, trading, and
  volatility}, volume~1.
\newblock CRC Press, 2017.

\bibitem[Shang and You(2018)]{shang2018distributionally}
Chao Shang and Fengqi You.
\newblock Distributionally robust optimization for planning and scheduling
  under uncertainty.
\newblock \emph{Computers \& Chemical Engineering}, 110:\penalty0 53--68, 2018.

\bibitem[Shao and Zavala(2019)]{shao2019space}
Yue Shao and Victor~M Zavala.
\newblock Space-time dynamics of electricity markets incentivize technology
  decentralization.
\newblock \emph{Computers \& Chemical Engineering}, 127:\penalty0 31--40, 2019.

\bibitem[Street et~al.(2009)Street, Barroso, Granville, and
  Pereira]{street2009bidding}
Alexandre Street, Luiz~Augusto Barroso, Sergio Granville, and Mario~Veiga
  Pereira.
\newblock Bidding strategy under uncertainty for risk-averse generator
  companies in a long-term forward contract auction.
\newblock In \emph{2009 IEEE Power \& Energy Society General Meeting}, pages
  1--8. IEEE, 2009.

\bibitem[Wang et~al.(2016)Wang, Glynn, and Ye]{wang2016likelihood}
Zizhuo Wang, Peter~W Glynn, and Yinyu Ye.
\newblock Likelihood robust optimization for data-driven problems.
\newblock \emph{Computational Management Science}, 13\penalty0 (2):\penalty0
  241--261, 2016.

\bibitem[Xavier et~al.(2021)Xavier, Qiu, and Ahmed]{xavier2021learning}
{\'A}linson~S Xavier, Feng Qiu, and Shabbir Ahmed.
\newblock Learning to solve large-scale security-constrained unit commitment
  problems.
\newblock \emph{INFORMS Journal on Computing}, 33\penalty0 (2):\penalty0
  739--756, 2021.

\bibitem[Xie(2020)]{xie2020tractable}
Weijun Xie.
\newblock Tractable reformulations of two-stage distributionally robust linear
  programs over the type-$\infty$ {W}asserstein ball.
\newblock \emph{Operations Research Letters}, 48\penalty0 (4):\penalty0
  513--523, 2020.

\bibitem[Yau et~al.(2011)Yau, Kwon, Rogers, and Wu]{yau2011financial}
Sheena Yau, Roy~H Kwon, J~Scott Rogers, and Desheng Wu.
\newblock Financial and operational decisions in the electricity sector:
  Contract portfolio optimization with the conditional value-at-risk criterion.
\newblock \emph{International Journal of Production Economics}, 134\penalty0
  (1):\penalty0 67--77, 2011.

\bibitem[Zavala(2013)]{zavala2013real}
Victor~M Zavala.
\newblock Real-time optimization strategies for building systems.
\newblock \emph{Industrial \& Engineering Chemistry Research}, 52\penalty0
  (9):\penalty0 3137--3150, 2013.

\bibitem[Zhang and Grossmann(2016)]{zhang2016enterprise}
Qi~Zhang and Ignacio~E Grossmann.
\newblock Enterprise-wide optimization for industrial demand side management:
  Fundamentals, advances, and perspectives.
\newblock \emph{Chemical Engineering Research and Design}, 116:\penalty0
  114--131, 2016.

\end{thebibliography}

\end{document}